\RequirePackage{fix-cm}
\documentclass[smallextended,numbook,runningheads]{svjour3}
\smartqed  
\usepackage{color,xcolor}
\usepackage{amsmath,amssymb,bm}
\usepackage{graphicx,graphics}
\definecolor{refgreen}{rgb}{0,0.5,0}
\definecolor{refgreen2}{rgb}{0.2,0.7,0}
\usepackage[colorlinks,linkcolor=blue,citecolor=blue]{hyperref}
\usepackage{booktabs,subfig} 
\usepackage{dsfont}
\usepackage{yfonts}
\usepackage{array}
\usepackage{makecell} 

\usepackage{tikz}
\usepackage{pgfplots}
\usetikzlibrary{arrows}
\usetikzlibrary{decorations.pathreplacing}
\usepackage{adjustbox}
\usetikzlibrary{intersections}
\usepackage{leftidx}
\usetikzlibrary{calc}

\usepackage{amsthm}
\newtheorem{thm}{Theorem}[section]
\newtheorem{lem}[thm]{Lemma}

\newtheorem{rem}[thm]{Remark}

\newcommand\calH{{\cal H}}
\newcommand\calG{{\cal G}}

\newcommand\calT{{\mathcal T}}

\newcommand\calE{{\mathcal E}}
\newcommand\calV{{\mathcal V}}

\newcommand\calI{{\mathcal I}}
\def            \d          {\hspace{1.5pt}\mathrm{d}}
\newcommand{\D}{{\textnormal d}}
\newcommand{\R}{\mathbb{R}}

\newcommand{\ddt}{\frac{\D}{\D t}}

\newcommand{\Ga}{\varGamma}
\newcommand{\Om}{\varOmega}
\renewcommand{\nu}{\textnormal{n}}

\newcommand{\comesh}{{n-1 \oplus n}} 

\newcommand\eps{\epsilon}
\newcommand\inveps{\eps^{-1}}

\newcommand\TOL{\textnormal{TOL}}

\newcommand\forquad{\quad\hbox{ for }\quad}
\newcommand\forallquad{\quad\hbox{ for all }\quad}

\journalname{...}
\date{ \phantom{b} \vspace{45mm}\phantom{e}}
\def\makeheadbox{\hspace{-4mm}Version of \today}

\begin{document}

\title{A posteriori error estimates for the bulk--surface Cahn--Hilliard system with dynamic boundary conditions}

\titlerunning{A posteriori error estimates for the bulk--surface Cahn--Hilliard system}        

\author{Nils Bullerjahn, Michael Lantelme}

\authorrunning{N.Bullerjahn and M.Lantelme} 

\institute{
	Nils Bullerjahn \at
	Institute of Mathematics, Paderborn University, \\
	Warburgerstr.~100, 33098 Paderborn, Germany, \email{bullerja@math.uni-paderborn.de} \\
	\and
	Michael Lantelme \at
	Institute of Mathematics, Paderborn University, \\
	Warburgerstr.~100, 33098 Paderborn, Germany, \email{lantelme@math.uni-paderborn.de} 
}

\date{\today}

\maketitle

\begin{abstract}
	This paper develops reliable residual-based a posteriori error estimators for the full discretization of a bulk--surface Cahn--Hilliard system with a two-parameter family of dynamic boundary conditions on a curved domain. The main novelty are theoretically justified computable temporal, coarsening and geometric error estimators. The numerical method combines a linear bulk--surface finite element method in space and a linearly implicit backward Euler method in time. Based on the derived error indicators, a space--time adaptive algorithm is proposed and numerical experiments illustrate and complement the theoretical result.
\end{abstract}

\newpage
\section{Introduction}

In this paper we develop an residual-based a posteriori error analysis for the bulk--surface Cahn--Hilliard system with a two-parameter family of dynamic boundary conditions in a curved domain with numerically computable error estimators, which are tested in a classical space--time adaptive algorithm. The numerical discretization combines a non-conforming linear bulk--surface finite element discretization in space, following the definition in \cite{ElliottRanner2013}, and a linearly implicit backward Euler method in time. Optimal-order a priori error estimates for the uniform numerical scheme are available in \cite{Bullerjahn2025}. 

The main novel contribution of this paper is the adaptation of the geometric error estimation techniques developed for surfaces in \cite{KovacsLantelme2026}, respectively evolving surfaces in \cite{Lantelme2026}, to the case of a curved flat domain. This is achieved by carefully combining the arguments on surfaces with the construction and estimates for the linear bulk--surface finite element method developed in \cite{ElliottRanner2013}. As a model problem we have chosen a bulk--surface Cahn--Hilliard system with dynamic boundary conditions. We believe that for this problem with dynamic boundary conditions it is especially interesting to be able to deal with a curved domain, as the interaction with the boundary is of particular interest. Furthermore, we are only aware of one work dealing with a posteriori error estimates for a linear parabolic model problem \emph{with dynamic boundary conditions}, dealing with a particular reformulation as a partial-differential--algebraic equation and only focusing on refinement, see \cite{AltmannZimmer2024}. So this work also serves as a first example of how to deal with dynamic boundary conditions in an a posterior error analysis with a standard refinement and coarsening strategy.

\smallskip
The classic Cahn-Hilliard equation was introduced in \cite{CahnHilliard1958} as a diffuse interface model for the phase separation process of binary mixtures. However, the classic Neumann boundary conditions lead to some practical limitations of the model, for example the contact angle between the interface between the two mixtures and the boundary is forced to $\pi/2$, see, e.g., \cite{Jacqmin2000}. In order to overcome these limitations several dynamic boundary conditions were proposed. In this paper we focus on the recently introduced two-parameter family of Cahn--Hilliard-type dynamic boundary conditions in \cite{KnopfStange2024}, containing most of the previously developed dynamic boundary conditions of Cahn--Hilliard-type \cite{GoldsteinMiranvilleSchimperna2011,LiuWu2019,KnopfLam2020,KnopfLamLiuMetzger2021} as special cases, cf. Figure~\ref{fig:graphDynBC}.

\smallskip
The general framework of obtaining a posteriori error estimates for parabolic problems was developed in \cite{Verfuerth2003}. Since then, several works have dealt with the reliability analysis for Cahn--Hilliard equations, analyzing semi-discretizations on polyhedral domains: In \cite{FengWu2007} the authors focused on reducing the classically exponential order dependence of the constants on the diffuse interface parameter $\eps$ to a polynomial order by using spectral estimates when no topological changes appear in the interface between the two phases. This was later extended to more general Cahn--Hilliard equations in \cite{BartelsMuller2010} and allowing for some topological changes in the interface in \cite{BartelsMullerOrtner2011}. Results for the more challenging obstacle potential are obtained in \cite{BanasNurnberg2009} and \cite{BartelsMuller2011}. To our knowledge this is the first fully discrete proof of reliability for a Cahn--Hilliard equation on a curved domain, introducing novel theoretically justified temporal, coarsening and geometric estimators, along with the classical spatial estimators adapted to incorporate the dynamic boundary conditions. We illustrate and test the efficiency of our theoretical results by a space--time adaptive algorithm. As far as we know, theoretical efficiency results for Cahn-Hilliard problems are still an open problem, see \cite[Remark~4.2]{BanasNurnberg2009}.

\smallskip
The paper is structured as follows: In Section~\ref{section:Cahn--Hilliard intro} we introduce the bulk--surface Cahn--Hilliard system with dynamic boundary conditions, which is the starting point of our analysis, and in Section~\ref{sec:discretization} we introduce the full discretization in space by bulk--surface finite elements and in time by the backward Euler method and obtain the core new geometric estimate. Section~\ref{sec:Reliability} then presents the main reliability error analysis, while in Section~\ref{sec:Numerical Experiments} we report on a number of numerical experiments to illustrate and complement our theoretical results.

\section{Bulk--surface Cahn--Hilliard system with dynamic boundary conditions \label{section:Cahn--Hilliard intro}}

This section is devoted to the introduction of the Cahn--Hilliard system with a two-parameter family of Cahn--Hilliard-type dynamic boundary conditions which is the starting point of our numerical analysis. The fully general system, even with additional convection terms, was first introduced in \cite{KnopfStange2024} and generalizes the important previously developed dynamic boundary conditions in \cite{GoldsteinMiranvilleSchimperna2011}, \cite{LiuWu2019}, \cite{KnopfLam2020}, and \cite{KnopfLamLiuMetzger2021}; cf. Figure~\ref{fig:graphDynBC}.

\subsection{Bulk--surface Cahn--Hilliard system with dynamic boundary conditions and asymptotic limits}
Let $\Om \subset \R^d$, with $d\in\{2,3\}$, be a bounded domain with a $C^2$-boundary $\Ga =\partial \Om$ and $T >0$ a final time. Then we introduce the bulk--surface phase field 
$$u=(\phi,\psi) \colon  (\Om \times [0,T]) \times(\Ga \times [0,T]) \to \R$$
and the bulk--surface chemical potential $$w=(\mu,\theta) \colon (\Om \times [0,T]) \times(\Ga \times [0,T]) \to \R.$$
Then the bulk--surface Cahn--Hilliard system with dynamic boundary conditions is given by the classical Cahn--Hilliard equation for the bulk phase field $\phi$ and bulk chemical potential $\mu$, c.f. \cite{CahnHilliard1958}, coupled to analogous Cahn--Hilliard-type dynamics for the surface phase field $\psi$ and surface chemical potential $\theta$:
\begin{subequations}\label{eq:CHreact}
	\begin{alignat}{2}
		\partial_t \phi =&\  m_\Om \Delta \mu \qquad&&\text{ in } \Om \times (0,T), \\
		\mu =&\  - \eps \Delta \phi + \eps^{-1} F_\Om'(\phi) &&\text{ in } \Om \times (0,T), \\
		\partial_t \psi =&\  m_\Ga \Delta_\Ga \theta - \beta m_\Om \partial_\nu \mu  &&\text{ on } \Ga \times (0,T), \\
		\theta =&\  - \delta \kappa \Delta_\Ga \psi + \delta^{-1} F_\Ga'(\psi) + \alpha \eps \partial_\nu \phi \qquad &&\text{ on } \Ga \times (0,T),
	\end{alignat}
for the initial value $u_0=(\phi_0,\psi_0) \colon \Om \times \Ga \to \R$ and where $m_\Om,m_\Ga >0$ denote the constant diffusion mobility, $\eps,\delta,\kappa >0$ denote diffusion interface parameters, and $F_\Ga,F_\Ga \colon \R \to \R$ are free energy potentials (e.g. a double well with minima $\pm 1$). The differential operators on the surface are the standard Laplace--Beltrami operator $\Delta_\Ga$ and normal derivative $\partial_\nu$, as defined, e.g., in \cite{DziukElliott2013}. In this work we consider, for the two parameters $K,L \in (0,\infty)$, the general family of Robin-type coupling conditions,

\begin{align}
\eps K\partial_\nu \phi =&\  \alpha \psi - u, &&\text{ on } \Ga \times (0,T),\label{eq:robintypecond1}\\
L \partial_\nu \mu =&\ \beta \theta - \mu, &&\text{ on } \Ga \times (0,T), \label{eq:robintypecond2}.
\end{align}
\end{subequations}
where $\alpha, \beta \in \R$. 

The parameter $L$ in \eqref{eq:robintypecond2} is given as the proportionality constant between the mass flux, i.e. the motion of material towards and away from the boundary, and the difference of the bulk and surface chemical potentials, $1/L$ can be interpreted as a reaction rate. Similarly, $1/K$ can be interpreted as a transmission rate, corresponding, for instance, to a case where the bulk materials are transformed at the boundary by a chemical reaction.

\begin{figure}[hb] 
	\centering
	\begin{adjustbox}{width=0.8\textwidth}
		\begin{tikzpicture}
			\draw[-stealth, very thick] (0,0)  -- (5.4,0) node[right] {$K$}; 
			\node at (0,-0.3) {$0$};
			\node at (5,-0.3) {$\infty$};
			\draw[-stealth, very thick] (0,0)  -- (0,5.4) node[above] {$L$}; 
			\node at (-0.3,0) {$0$};
			\node at (-0.3,5) {$\infty$};
			
			\draw[dashed, very thick] (5,0) -- (5,5);
			\draw[dashed, very thick] (0,5) -- (5,5);
			
			\node[draw,align=center] at (-1,5.6) {LW-model \\ \cite{LiuWu2019}};
			
			\node[draw,align=center] at (2.5,5.6) {transmission rate dependent \\ dyn. bound. cond. \cite{KnopfLam2020}};
			
			\node[draw,align=center] at (-1.1,-0.6) {GMS-model \\ \cite{GoldsteinMiranvilleSchimperna2011}};
			
			\node[draw,rotate=90, align=center] at (-0.6,2.5) {reaction rate dependent \\ dyn. bound. cond. \cite{KnopfLamLiuMetzger2021}};
			
			\node[draw,align=center] at (6.6,5.6) {Neumann boundary \\ conditions \cite{CahnHilliard1958}};

			\fill[gray!20] (0.01,0.01) rectangle (4.99,4.99);
			
			
			\node[draw,align=center] at (2.5,2.5) {bulk--surface system with \\ dyn. bound. cond. \cite{KnopfStange2024}};
			
			
		\end{tikzpicture}
	\end{adjustbox}
	\caption{Sketch of the different types of dynamic boundary conditions corresponding to the parameter-choices $K,L \in [0,\infty]$. Adapted from \cite{Bullerjahn2025}.}\label{fig:graphDynBC}
\end{figure}
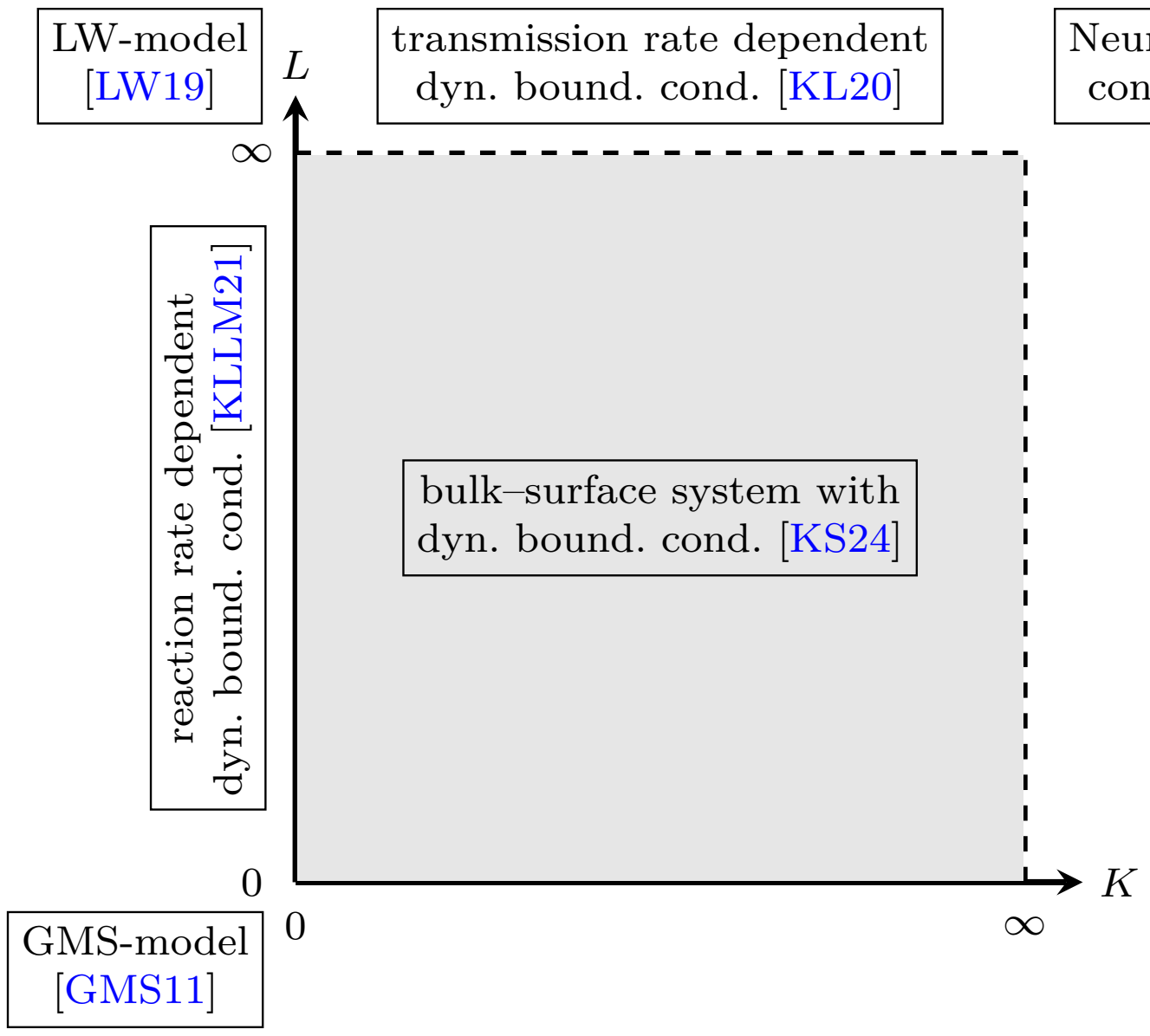

Note that the Robin-type coupling conditions \eqref{eq:robintypecond1} \& \eqref{eq:robintypecond2} also serve as an interpolation between Dirichlet- and Neumann-type boundary conditions for the phase fields and chemical potentials. For example in the limit case of instantaneous reaction, i.e. $L=0$, we recover the Dirichlet boundary condition in the chemical potentials and on the other hand a vanishing reaction rate, i.e. $L=\infty$, results in the Neumann boundary condition, and similar for $K$. So this system combines several previously studied coupling conditions, see \cite{GoldsteinMiranvilleSchimperna2011}, \cite{LiuWu2019}, \cite{KnopfLam2020}, and \cite{KnopfLamLiuMetzger2021}, for different choices of $K$ and $L$ as illustrated in Figure~\ref{fig:graphDynBC}. It has also been shown on the level of weak solutions in \cite{KnopfStange2024}, that for $K$ and $L$ converging to either $0$ or $\infty$, solutions of the system \eqref{eq:CHreact} converge to solutions of the respective limit system.

\begin{rem} \label{rem:Cont results Epsilon}
	It has been shown in \cite{Pego1989}, that for the classical Cahn--Hilliard equation (i.e. the choice $K=L=\infty$) and a typical double well potential with minima at $\pm 1$, the diffusion interface parameter $\eps >0$ governs the thickness of internal transition layers and that for $\eps \ll 1$ the phase field attains values close to $\pm 1$ corresponding to the pure phases of the materials in a time of $O(\eps^2)$. We strongly expect this typical behaviour to also transfer to the present model and the diffusion interface parameters $\eps, \delta, \kappa >0$, as it has been observed in numerous numerical examples across the literature.
	
	After this initial phase it has been shown that the solution evolves on the dominant slowest time scale to a steady state with an interface velocity determined by the mean curvature of the interface, see \cite{Pego1989} and \cite{KnopfPoiattiStangeYayla2026}.
\end{rem}

As we are interested in developing a time and space adaptive algorithm for this problem, Remark~\ref{rem:Cont results Epsilon} suggest that the diffusion interface parameters play an important role in the error indicators, since we expect spatial refinements especially around the internal transition layers. And in fact, it has been investigated in previous works for similar Cahn--Hilliard problems, see e.g. \cite{FengWu2007}, that the choice of the diffusion mobility parameters $m_\Om$ and $m_\Ga$ and especially the diffusion interface parameters $\eps, \delta$ and $\kappa$ has a substantial effect on the constants appearing in the a posteriori analysis. To investigate how the resulting estimators depend on these parameters, we focus our analysis on a single parameter, $0< \eps<1$, by the choice 
\begin{equation}
	m_\Om=m_\Ga=\kappa=1 \text{ and } 0<\eps=\delta<1, \label{eq:Choice of eps delta...}
\end{equation}
while capturing the most important relations of the general case.

\subsection{Abstract weak formulation}
For $1 \leq p \leq \infty$ and $k\geq 0$ we denote the standard Lebesgue spaces on $\Om \subset \R^d$ by $L^p(\Om)$ and the Sobolev spaces by $W^{k,p}(\Om)$ (see, e.g., \cite[Section~5.2.2]{Evans1998}). For $p=2$ we use the standard notation for the Hilbert spaces $H^k(\Om)=W^{k,2}(\Om)$. For the Lebesgue and Sobolev spaces on $\Ga=\partial \Om$ we use the similar notation $L^p(\Ga)$ and $W^{k,p}(\Ga)$ (see, e.g., \cite[Definition 2.11]{DziukElliott2013}). For any Banach Space $X$, we denote by $\| \cdot \|_X$ its norm, by $X'$ the dual space and the associated duality pairing by $\langle \cdot , \cdot \rangle_X$. We introduce the combined Hilbert spaces
\begin{align}
	&\calH:= L^2(\Om) \times L^2(\Ga), \quad  \calV:= H^1(\Om) \times H^1(\Ga), \label{eq:Def analytic function spaces}
\end{align}
with the standard combined scalar products. Following the notation from \cite[Section~2.1.1]{KovacsLubich2017}, we introduce the following bilinear forms and semi-norm
\begin{equation} \label{eq:Def bilinear forms}
	\begin{aligned}
		&m \colon \calH \times \calH \to \mathbb{R} ,\\
		&\quad ((\phi,\psi),(\zeta,\xi)) \mapsto  \int_\Om \phi  \zeta + \int_\Ga  \psi \xi  , \\
		&a^{J,\lambda,\eps} \colon \calV \times \calV \to \mathbb{R} ,\\
		&\quad ((\phi,\psi),(\zeta,\xi)) \mapsto  \eps \int_\Om \nabla \phi \cdot \nabla \zeta + \eps \int_\Ga \nabla_\Ga \psi \cdot \nabla_\Ga \xi  + J^{-1} \int_\Ga (\lambda \psi-\phi)(\lambda \xi-\zeta) ,\\
		&|\cdot|_{a^{J,\lambda,\eps}} \colon \calV \to \mathbb{R} ,\\
		&\quad (\phi,\psi) \mapsto  \Big(a^{J,\lambda,\eps}((\phi,\psi),(\phi,\psi))\Big)^{\frac12}. 
	\end{aligned}
\end{equation}

We additionally define the weak dual norms, for $d \in \calH$:
\begin{equation} \label{eq:dualNormDef} 
	\begin{aligned}
		\| d \|_{\ast}:=\| d \|_{\calV^{\ast}}:= \sup_{0\not=\varphi \in \calV}\frac{m(d,\varphi)}{\|\varphi\|_{\calV}}.
	\end{aligned}
\end{equation}

Following \cite{KnopfStange2024}, a weak solution $(u,w)$ of the system \eqref{eq:CHreact} on $[0,T]$, for $K,L \in (0,\infty)$ and $u_0 \in V$, satisfies the following properties:
\begin{enumerate}
	\item[i)] The functions $u$ and $w$ have the regularity:
	\begin{equation} \label{eq:WeakFormReg}
		\begin{aligned}
			&u \in C([0,T];\calH) \cap H^1(0,T;(\calV)') \cap L^\infty(0,T;\calV), \\
			&w \in L^2(0,T;\calV). 
		\end{aligned} 
	\end{equation} 
	\item[ii)] The function $u$ satisfies the initial conditions $u|_{t=0}=u_0$ a.e. in $\Om \times \Ga$.
	\item[iii)] The functions $u$ and $w$ satisfy the weak formulation
	\begin{subequations} \label{eq:CHweak}
		\begin{align}
			m(\partial_t u,\zeta) + a^{L,\beta}(w,\zeta) =&\ 0, \label{eq:CHweak1}\\
			m(w,\eta) - a^{K,\alpha,\eps}(u,\eta) =&\  m(\inveps F'(u),\eta),\label{eqCHweak2}
		\end{align}
	\end{subequations}
	a.e. in $[0,T]$, for all $\eta,\zeta \in V$, where we abbreviated $F'(u) = (F_\Om'(\phi),F_\Ga'(\psi))$, for $u=(\phi,\psi)$ and denoted $a^{L,\beta}:=a^{L,\beta,1}$.
\end{enumerate}

The main result in \cite[Theorem 3.2]{KnopfStange2024} provides existence of weak solutions under polynomial growth conditions on the nonlinear potentials, and even a higher regularity result.

\begin{rem}
	In the following, we mostly use a bulk--surface notation, e.g. $u$ instead of $\phi$ and $\psi$, if no confusion or difficulties can arise, in order to keep the formulas shorter. Note however, that there are some subtle differences appearing between bulk and surface were we switch back to the longer notation.
\end{rem}

\section{Full discretization by bulk--surface finite elements and the backward Euler method}
\label{sec:discretization}
In order to discretize the equation in space, we use the linear bulk--surface finite element method for partial differential equations developed by Elliott and Ranner \cite[Sections 4,5]{ElliottRanner2013}, see also \cite[Section~4]{KovacsLubich2017} and \cite[Sections 4-7]{ElliottRanner2021}. The discretization in time is done via the backward Euler method. \\

\textbf{Spatial discretization}
Approximate the domain $\Om$ by a polyhedral domain $\Om_h$ and take a shape-regular triangulation $\calT_h$ of $\Om_h$ (see \cite[Definition 4.4.13]{BrennerScott2008}). This triangulation consists of closed simplices, either triangles in $\mathbb{R}^2$ or tetrahedra in $\mathbb{R}^3$, and we denote by $E_h$ the set of $(d-1)$-dimensional sides of elements in $\calT_h$. Let $\Ga_h:=\partial \Om_h$, then the faces of $\Ga_h$ are $(d-1)$-simplices, hence edges in $\mathbb{R}^2$ or triangles in $\mathbb{R}^3$, and we assume that the vertices lie on $\Ga$, and we denote by $E_{\Ga_h}$ the set of boundary edges and by $\mathcal{P}_h$ the set of $(d-2)$-dimensional sides of the faces in $\Ga_h$. We define the local mesh width $h_T:=\text{diam}(T)$ for $T \in \calT_h$ and the mesh width as $h=\max\{\text{diam}(T) \colon T \in \calT_h\}$ and assume that the mesh width is sufficiently small, such that for all $x \in \Ga_h$, there exists a unique point $p(x) \in \Ga$, such that $x-p(x)$ is orthogonal to the tangent space of $\Ga$ at $p$, for more details on the construction see \cite[Section~4]{ElliottRanner2013}. We additionally define the local mesh width on $S \in E_h$ as $h_S:=\text{diam}(S)$, and on $P \in \mathcal{P}_h$ as $h_P:=\text{diam}(P)$, note that for $d=2$ we have $h_P=0$.

We define the non-conforming finite element space $V_h\nsubseteq H^1(\Om)$, corresponding to $\calT_h$, by the continuous, piecewise linear nodal basis functions on $\Om_h$. Similarly the nodes of the triangulation $\calT_h$ on the boundary $\Ga_h$ define the \emph{surface} finite element space $S_h\nsubseteq H^1(\Ga)$. Refer to \cite[Section~6.2]{ElliottRanner2021} for more details.

The discrete gradient $\nabla_h$ and tangential gradient $\nabla_{\Ga_h}$ are piecewisely defined analogously to $\nabla$ and $\nabla_\Ga$, see \cite[Section~4]{DziukElliott2013}. The discrete trace operator $\gamma_h v_h$ is defined by the restriction of the continuous functions $v_h \in V_h$ onto $\Ga_h$, with the same notational conventions as before, e.g. $\nabla_{\Ga_h}v_h:=\nabla_{\Ga_h}(\gamma_h v_h)$.

For the combined bulk--surface system we define analogous to the analytic spaces in \eqref{eq:Def analytic function spaces} the \emph{bulk--surface} finite element space as
\begin{align*}
	\calV_{h}:= V_h \times S_h.	
\end{align*}

\begin{rem}
\label{remark:NVB}
To be able to define and analyse the temporal discretization and the mesh-transfer operator, we specify the refinement and coarsening strategy.
We highlight, that this will be particularly required for the statement and proof of Lemma~\ref{Lemma:common_interpolation}.
In accordance to the assumption in \cite{KovacsLantelme2026} we will assume that any mesh is generated based on some initial admissible coarse mesh such that the lift on the boundary is a bijective map. 
Additionally the refinement is done in the following manner: After each refinement, nodes on the boundary are immediately lifted to the boundary; at most one node for each triangle has to be lifted; and coarsening only effects the deepest refinement level. Note that this does not restrict multiple refinement or coarsening of an element, but dictates the order in which the mesh changes take effect.
This is an essential feature as it helps to guarantee that hierarchically identical vertices always align between time steps.
Throughout the paper we assume the use of Newest-Vertex-Bisection (NVB), see \cite{Mitchell1991}, combined with immediate lifting, which fulfills all requirements.
\end{rem}

\textbf{Temporal discretization}
We divide the interval $[0,T]$ into $N$ sub-intervals $[t^{n-1},t^n]$ of length $\tau_n$ for $n = 1,\dots,N$, where $0=t^0<\cdots <t^N = T$. In each time step we assume that a conforming triangulation $\Om_h^n$ is given, as described above. Notice here that even though the domain $\Om$ is stationary, the approximated domains $\Om_h^n$ can change, since after each refinement we lift the vertices to the boundary of the exact surface $\Ga$. We denote the corresponding triangulation by $\calT_h^n$ and the bulk--surface finite element spaces by $\calV_h^n:=V_h^n \times S_h^n$.

In order to compare finite element functions on the domain $\Om_h^{n-1}$ to functions on $\Om_h^{n}$ we introduce the mesh-transfer operator $$\calI_{n-1}^n\colon \calV_h^{n-1} \to \calV_h^n.$$
This interpolation operator is defined by the following idea: When we refine a triangle in $\Om_h^{n-1}$, the new node is initially within the original domain $\Om_h^{n-1}$, where we can evaluate a nodal value based on the finite element function. If the new node is on $\Ga_h^{n-1}$ we lift it to $\Ga$. The mesh-transfer operator is then defined by assigning these new nodal values with the prior evaluation in $\Om_h^{n-1}$. When we coarsen a triangle, a vertex is removed and its nodal value is discarded. See \cite{KovacsLantelme2026} who analysed the same situation on surfaces.  

Then we define the following backward finite difference for $\big(u_h^{n}\in \calV_h^{n}\big)_{n=0}^N$ as,
\begin{align}
	\partial_{\tau_n} u_h^n =\frac1{\tau_n} \big(u_h^n - \calI_\Om^n u_h^{n-1}\big) \in \calV_h^n.
\end{align}
We highlight, that the construction of the mesh-transfer interpolant yields a backwards difference which yields functions within the finite element space, which simplifies computation of the full discretization heavily and avoids the need to work with two distinct discrete domains.

\begin{rem} \label{rem:mass conservation}
Note that the mesh-transfer operator in the presented form does not ensure mass conservation.
Whereas the effect of discarding nodal values clearly effects the mass, there is a second source of mass change related to the lift.
On the boundary, due to the potential movement of vertices the elements are perturbed to better resolve the geometry, which introduces a geometric, high-order, local change in mass. 

The nodal definition allows us to simply generalize results required to quantify the interpolation errors in Lemma~\ref{Lemma:common_interpolation}, however extension to mass conserving interpolation operators could be of interest.
\end{rem}

\textbf{Full discretization}
Then we use a Galerkin ansatz in order to solve this variational problem in a finite element space and combine this with the discretization in time by the backward Euler method. The fully discrete scheme reads: Let $L,K \in (0,\infty)$, then we want to find for $n = 1,\dots, N$, the family of functions $u_h^n \in \calV_h^{n}$ and $w_h^n \in \calV_h^{n}$ such that
\begin{subequations} \label{eq:CHweakD}
	\begin{align}
		m_h^n(\partial_{\tau_n} u_h^n,\zeta_h) + a_h^{n,L,\beta}(w_h^n,\zeta_h) =&\ 0, \label{eq:CHweakD1}\\
		m_h^n(w_h^n,\eta_h) - a_h^{n,K,\alpha,\eps}(u_h^n,\eta_h) =&\  m_h^n(\inveps \calI_{n-1}^n F'( u_h^{n-1}),\eta_h).\label{eq:CHweakD2}
	\end{align}
\end{subequations}
for all $\zeta_h \in \calV_h^{n}$, $\eta_h \in \calV_h^{n}$ and with some given initial value $u_h^0 \in \calV_h^{0}$. 

In the following we will neglect the time dependence of the norms and bilinear forms in the notation, whenever it is clear from the context.

\subsection{Lift to the exact domain \& comparison with the exact solution}
\label{sec:lift_interpolation}
We require two steps to compare the set of discrete solutions obtained in \eqref{eq:CHweakD1}, \eqref{eq:CHweakD2} to the exact solution of \eqref{eq:CHweak1}, \eqref{eqCHweak2}.
Namely these are lifting, i.e., mapping functions from the discrete to the continuous domain, and time interpolation, to extend the fully discrete sequence of solutions on the full time-scale.
The upcoming sections hold for dimension $d \geq 2$.

\textbf{Lifting} 
We introduce the lift operator as described in \cite[Section~4]{ElliottRanner2013}, following the construction of Dubois \cite{Dubois1990} and notation of \cite{Deckelnick2009}. 
	We define the following objects:
	\begin{center}
		\begin{tabular}{  m{2cm}  m{9.2cm} } 
			$T$&  element, i.e. $d$-dimensional simplex in $\calT_h^n$\\ 
			$\hat T$ & reference-simplex in $\R^d$\\ 
			$T^\ell$ & curved element corresponding to $T$, such that $(\calT_h^n)^\ell =\Om$\\ 
			$F_T \colon \hat T \to T$ & affine reference map $F_T(\hat x) = A \hat x + b$ for $\hat x \in \hat T$ \\
			$F_T^e \colon \hat T \to T^\ell$ & reference map of a curved simplex, given by $F_t^e = F_T + \Phi$ \\
			$G_h \colon T \to T^\ell$ &  lift map of the element $T$ given by $F_T^e \circ (F_T)^{-1}$ \\
			$x_1,\dots,x_{d+1}$ &  nodes of an element $T$, where the, potential empty, set of boundary nodes are always labeled first, i.e. $x_i \in \Ga_h^n$ for $i = 1,\dots,l$, where $l$ denotes the number of boundary nodes.
		\end{tabular}
	\end{center}

First note that for boundary points $x \in \Gamma_h$ we employ the standard lift operator, i.e. the closest point projection $p(x) := x^\ell$, introduced in \cite{Dziuk1988}. 
Secondly for elements in the bulk $T \subset \Omega_h$ only those with at least one boundary edge are affected by the lift. Represent $x \in T$ by its barycentric coordinates $x = \sum_{i=1}^{d+1} \hat{\lambda}_i x_i$ and define $\hat{\lambda}^\star = \sum_{i=1}^{l} \hat{\lambda}_i$.
	This allows us to define the projection from an interior point of the simplex to its boundary edge in $\Gamma_h$ by $y(\hat{x}) = \sum_{i=1}^l \frac{\hat{\lambda}_i}{\hat{\lambda}^\star} x_i$.
Now as for the boundary point case we utilize the normal projection $p(y) \in \Gamma$ to define the map $$\Phi(\hat{x}) := (\hat{\lambda}^\star)^{3} (p(y)-y),$$ for non-fixed points of the bulk. 
The bulk lift is then summarized as the map $G_h$ as a combination of the affine reference map $F$ and the non-linear map $\Phi$.
	
The process of lifting functions is understood for any function $\eta_h \in V_h$ as 
	\begin{equation*}
		\eta_h^{\ell} := \eta_h \circ G_h^{-1}.
	\end{equation*} 

We write the lifted finite element space as 
\begin{align*}
	\calV_h^\ell:=\{v_h^\ell : v_h \in \calV_h\}.
\end{align*}

Due to the change in domains the lift from different discrete time steps changes. We write the additional dependency by a subscript, e.g.\ for $\Om_h^n$ its corresponding lift is $(\cdot)^{\ell_n}$.

The following geometric approximation error bounds, coming from the arguments in the proof of \cite[Lemma~6.2]{ElliottRanner2013}, estimate the difference of the lifted and unlifted $L^2$- and $H^1$-inner products:
\begin{lem}[geometric approximation errors] \label{Lemma:GeomAppError}
	The following estimates hold for sufficiently small $h$ and for any $v,w \in H^1(\Om)$ and $\phi,\eta \in H^1(\Ga)$,
	\begin{align*}
		\Big|\int_\Om \nabla v \cdot \nabla w - \int_{\Om_h^n}\nabla_h v^{-\ell_n} \nabla_h w^{-\ell_n}\Big| \leq&\  C \sum_{T \in B_h^n} h_T \int_{T} \big|\nabla v^{-\ell_n} \cdot \nabla w^{-\ell_n} \big|, \\
		\Big|\int_\Ga \nabla_\Ga \phi\cdot \nabla_\Ga \eta - \int_{\Ga_h^n}\nabla_{\Ga_h^n}\phi^{-\ell_n}\cdot \nabla_{\Ga_h^n} \eta^{-\ell_n} \Big| \leq&\ C \sum_{S \in E_{\Ga_h}^n } h_S^2 \int_{S}\big|\nabla_{\Ga_h^n}\phi^{-\ell_n}\cdot \nabla_{\Ga_h^n} \eta^{-\ell_n} \big|, \\
		\Big|\int_\Om v w  - \int_{\Om_h^n} v^{-\ell_n} w^{-\ell_n} \Big| \leq&\   C \sum_{T \in B_h^n} h_T \int_{T} \big|v^{-\ell_n} w^{-\ell_n} \big|, \\
		\Big|\int_\Ga \phi \eta  - \int_{\Ga_h^n} \phi^{-\ell_n} \eta^{-\ell_n} \Big| \leq&\  C \sum_{S \in E_{\Ga_h}^n} h_S^2 \int_{S}\big|\phi^{-\ell_n} \eta^{-\ell_n} \big|,
	\end{align*}
	where $B_h^n$ collects elements which have at least two nodes on the boundary, and recall that $E_{\Ga_h}^n$ is the set of boundary edges.
\end{lem}

\begin{rem} \label{rem:normequivalence}
	As a consequence we also have the $h$-uniform equivalence of the norms $\| \cdot \|_{H^1(\Om)}$ and $\| \cdot \|_{H^1(\Om_h^n)}$ and of the norms $\| \cdot \|_{L^2(\Om)}$ and $\| \cdot \|_{L^2(\Om_h^n)}$.
	Respectively for the surface norms the equivalence holds for $\| \cdot \|_{H^1(\Ga)}$ and $\| \cdot \|_{H^1(\Ga_h^n)}$ and of the norms $\| \cdot \|_{L^2(\Ga)}$ and $\| \cdot \|_{L^2(\Ga_h^n)}$.
\end{rem}

\begin{figure}[hb] 
	\begin{adjustbox}{width=\textwidth}
		\begin{tikzpicture}
			
			\draw[-] (0,4) -- (6,4);
			\draw[-] (3,1) -- (6,4);
			\draw[-] (0,4) -- (3,1);
			
			\filldraw (0,4) circle (0.05); 
			\filldraw (6,4) circle (0.05);
			\filldraw (3,1) circle (0.05);
			\filldraw (1.5,4) circle (0.05);
			\filldraw (1.24,4.68) circle (0.05);
			\draw (-0.54,3.54) arc (135:45:5);
			
			\draw[dashed] (1.5,4) -- (1.24,4.68) node[midway,right] {$\ell_{n-1}$};
			\node[align=center,thick] at (1.5,3.7) {$x$};
			\node[align=center,thick] at (1.24,4.98) {$x_1$};
			
			\node[align=center, very thick] at (-0.8,3.5) {$\Ga$};
			
			\pgfmathsetmacro{\xshift}{8.5}
			\coordinate (A) at (0+\xshift,4);
			\coordinate (B) at (3+\xshift,5);
			\draw[-] (0+\xshift,4) -- (3+\xshift,5);
			\draw[-] (6+\xshift,4) -- (3+\xshift,5);
			\draw[-] (3+\xshift,1) -- (6+\xshift,4);
			\draw[-] (0+\xshift,4) -- (3+\xshift,1);
			\draw[-] (3+\xshift,5) -- (3+\xshift,1);

			\filldraw (0+\xshift,4) circle (0.05);
			\filldraw (6+\xshift,4) circle (0.05);
			\filldraw (3+\xshift,1) circle (0.05);
			\filldraw (3+\xshift,5) circle (0.05);
			\filldraw (1.5+\xshift,4.5) circle (0.05);
			\filldraw (1.42+\xshift,4.74) circle (0.05);
			\draw (-0.54+\xshift,3.54) arc (135:45:5);
			
			\node[align=center, text=gray] at (1.1+\xshift,4.88) {$x_1$};
			\filldraw[fill=gray!30, draw=gray!60] (1.24+\xshift,4.68) circle (0.05);
			
			\coordinate (M) at ($(A)!0.5!(B)$); 
			\fill (M) circle (0.05);
			\node[below left, xshift=0pt, yshift=0pt] at (M) {$\widetilde x$};
			\node[above right, xshift=0pt, yshift=0pt] at (M) {$\ell_{n}$};
			
			\draw[dashed] (M) -- (1.42+\xshift,4.74);
			\draw[dashed] (0+\xshift,4) -- (6+\xshift,4);
			\draw[dashed] (1.5+\xshift,4) -- (M); 
			\filldraw (1.5+\xshift,4) circle (0.05);
			\node[align=center,thick] at (1.5+\xshift,3.7) {$x$};
			\node[align=center,thick] at (1.42+\xshift,5.04) {$x_2$};
			
			\draw[->,thick] (\xshift-2,4) -- (\xshift-0.3,4) node[midway,above] {$\calI^{n}$};
			
			\node[align=center, very thick] at (-0.8+\xshift,3.5) {$\Ga$};

		\end{tikzpicture}
	\end{adjustbox} 
	\caption{Differences in lifting the domain from $\Om_{n-1}$ and $\Om_n$  to $\Om$ under to usage of the mesh-transfer interpolant $\calI_{n-1}^n$. Either we lift $x$ directly to $x_1$ via the lift $\ell_{n-1}$, or we correspond $x$ to $\widetilde x$ (i.e.\ fixing the midpoint and only lifting the vertices of the element) and then via $\ell_{n}$ to $x_2$, where in general, as in this example, $x_1 \not= x_2$.\label{Fig:ProblemWithLifts}}
\end{figure}
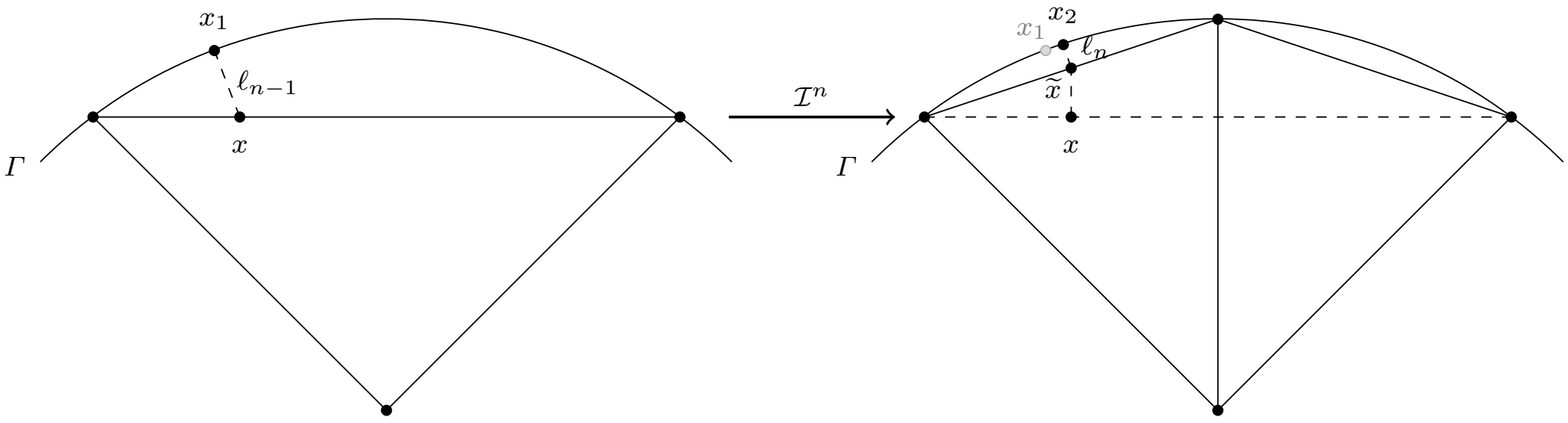

\textbf{Time interpolation}
We have to decide how we want to interpolate the discrete phase field values $u_h^n \in \calV_h^{n}$, obtained by the discrete scheme, to a function on the entire time interval. A first obvious choice is the following: linearly interpolate the lifted finite element functions to
\begin{align}
	u_{h,\tau}(t)=\frac{t-t^{n-1}}{\tau_n} (u_h^n)^{\ell_n} + \frac{t^n-t}{\tau_n} (u_h^{n-1})^{\ell_{n-1}}, \quad\text{ for } t \in (t^{n-1},t^n]. \label{eq:u_cont} 
\end{align}
This interpolation is continuous in time which is important for the continuous error estimates which we will perform, but has the disadvantage that it can not be unlifted to the finite element domain, since the lifts $\ell_n$ and $\ell_{n-1}$ are in general not identical as can be seen e.g.\ in Figure~\ref{Fig:ProblemWithLifts}. This would pose serious problems in the estimate of the residual by the error estimators.

In order to prevent this problem we propose a different interpolation: the piecewise linear interpolation in time is given by
\begin{align}
	\widetilde u_{h,\tau}(t)=\frac{t-t^{n-1}}{\tau_n} u_h^n + \frac{t^n-t}{\tau_n} \calI_{n-1}^n u_h^{n-1}, \quad\text{ for } t \in (t^{n-1},t^n]. \label{eq:u_computable}
\end{align}
This interpolation has the advantage that in each time interval $(t^{n-1},t^n]$, the interpolation $\widetilde u_{h,\tau}$ is contained in $\calV_h^{n}$. Unfortunately this interpolation is in general discontinuous in time at the discrete time levels, which means it is not suitable for the continuous error estimates. Hence, we will change to this interpolation in the estimation of the residual by the error estimators and have to estimate the difference of $\widetilde u_{h,\tau}^{\ell_n}$ and $u_{h,\tau}$, i.e., the difference of $(u_h^{n-1})^{\ell_{n-1}}$ and $(\calI_{n-1}^n u_h^{n-1})^{\ell_{n}}$.

Note that the two interpolations $(u_h^{n-1})^{\ell_{n-1}}$ and $(\calI_{n-1}^n u_h^{n-1})^{\ell_{n}}$ are in general not identical, their error can be decomposed into two parts: 
\begin{enumerate}
	\item By coarsening the mesh from $\calT_h^{n-1}$ to $\calT_h^{n}$, the information at the deleted grid points is lost in the interpolation $(\calI_{n-1}^n u_h^{n-1})^{\ell_{n}}$, but not in $(u_h^{n-1})^{\ell_{n-1}}$. This will lead to a coarsening error estimator in the a posteriori error estimates.
	\item The difference in the lifts $\ell_n$ and $\ell_{n-1}$ to the exact domain $\Om$ introduces a geometric error, even if there is no information lost due to the coarsening, as can be seen in Figure~\ref{Fig:ProblemWithLifts}.
\end{enumerate}

%
%
%

The situation for the chemical potential $w_h^n \in \calV_h^n$ is different. Since no derivative acts on this value, the discrete scheme is independent of the previous values of the chemical potential, hence we decided to interpolate this value piecewise constant:
\begin{align}
	w_{h,\tau}(t) = (w_h^n)^{\ell_n} , \forquad t \in (t^{n-1},t^n]. \label{eq:Def w_htau}
\end{align}
This choice is obviously not continuous, but since the continuous error estimates cancel this entirely, this poses no problems. Yet this choice is more suitable for the elliptic second equation \eqref{eq:CHweakD2}, as it eliminates the dependence of the estimators on the quantity $w_h^{n-1}$, which does not enter the numerical scheme in the current step.

\subsection{Smallest common triangulation for bulk--surface}
A common triangulation is essential to compare functions for adaptive meshes. Our setting will be based on the analysis of common triangulations of flat domains \cite[Section~2.2.2]{Kreuzer2013} and of surface domains in \cite{KovacsLantelme2026}. 

\textbf{Smallest common refinement.}
	Following \cite{KovacsLantelme2026} based on ideas of \cite[Section~2.2.2]{Kreuzer2013} we introduce the smallest common refinement of subsequent meshes $\Omega_h^{n-1}$ and $\Omega_h^n$, denoted by $\Omega_h^{\comesh} := \Omega_h^{n-1} \oplus \Omega_h^n$, which will serve as a common triangulation to compare functions on these meshes. 
	Recall the mesh assumptions stated in Section~\ref{sec:discretization}. Note that we assume that both subsequent meshes are in general not a refinement of one another and defined on different domains. But they are both a refinement of the initial mesh $\Omega_h^0$ obtained by NVB under possible hierarchical lifts, see Remark~\ref{remark:NVB}. In this sense the common mesh $\Omega_h^{\comesh}$ is also a refinement of $\Omega_h^0$.
	
	The smallest common refinement is defined by the following idea: $\Omega_h^{\comesh}$ is the collection of all common elements of the parent meshes, and always the refined part of $\Omega_h^n$ or $\Omega_h^{n-1}$ if an element was refined or coarsened.
	
\begin{figure}[htbp]
		\centering
		\includegraphics[scale=0.8]{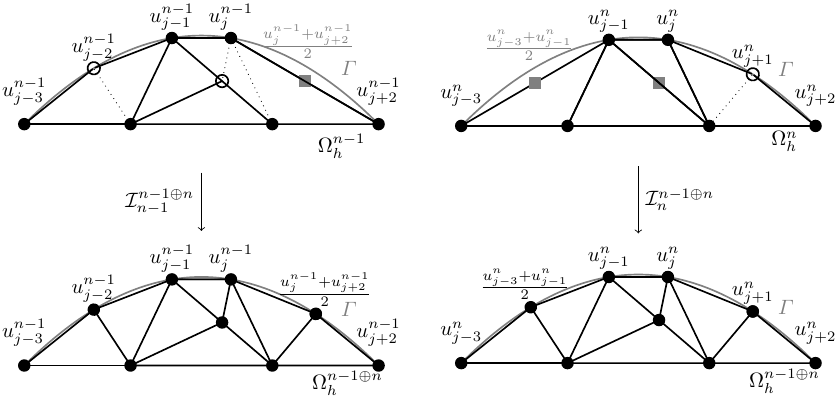}
		\caption{Construction of smallest common refinement $\Omega_h^{\comesh}$ and common interpolation operators $\mathcal{I}_k^{\comesh}$ ($k=n-1,n$). The $\bullet$ marks (top row) denote common nodes of $\Omega_h^{n-1}$ and $\Omega_h^n$. The common triangulation (bottom row) keeps the more refined parts of the parent triangulations. Where we highlighted refined edges by dotted lines. (The nodal values are numbered locally here)}
		\label{fig:nodal_values_of_common_triangulation}
	\end{figure}

	Nodes which appear in both parent meshes are referred to as common nodes of $\Omega_h^{n-1}$ and $\Omega_h^n$ (marked by $\bullet$ in the top row of Figure~\ref{fig:nodal_values_of_common_triangulation}). 
	
	Note that the interpolating points between common nodes, marked by \textcolor{gray}{\tiny $\blacksquare$} in the top row of Figure~\ref{fig:nodal_values_of_common_triangulation}, are not nodes of the mesh $\Omega_h^{n-1}$ or $\Omega_h^n$. We denote these meshes by $\widehat{\Omega}_h^{n-1}$ and $\widehat{\Omega}_h^n$, respectively, which by construction satisfy: the mesh obtained by lifting the nodes of either $\widehat{\Omega}_h^{n-1}$, or $\widehat{\Omega}_h^{n}$ yields the discrete surface $\Omega_h^{\comesh}$.  In fact, only the interpolating points needed to be lifted, i.e.~the ones marked by \textcolor{gray}{\tiny $\blacksquare$} in the top row of Figure~\ref{fig:nodal_values_of_common_triangulation}.
	
	Additionally, we will label the boundary of the common triangulation $\Omega_h^{\comesh}$ by $\Gamma_h^{\comesh}$.
	Its corresponding element set is, by our usual convention, given as $E_h^{\comesh}$.
	
	\textbf{Common finite element space.}
	The common finite element space $\calV_h^{\comesh}$  corresponding to the smallest common refinement $\Omega_h^{\comesh}$ is defined exactly as in Section~\ref{sec:discretization}.
	
	In order to map functions on $\Omega_h^{n-1}$ or $\Omega_h^n$ to functions on $\Omega_h^{\comesh}$, we define the common interpolation operators $\mathcal{I}_k^{\comesh} \colon \calV_{h}^{n} \rightarrow \calV_h^{\comesh}$ for both $k = n-1$ and $k=n$, respectively.
	
	The nodal values of $\mathcal{I}_k^{\comesh}$ are given by:
	\begin{enumerate}
		\item[--] All nodes which appear in both $\Omega_h^{n-1}$ and $\Omega_h^{\comesh}$ or $\Omega_h^n$ and $\Omega_h^{\comesh}$, respectively, keep their nodal values (e.g., nodes marked by $\bullet$ and $\circ$ in the top row of Figure~\ref{fig:nodal_values_of_common_triangulation}).
		\item[--] Nodes which are not present in either $\Omega_h^{n-1}$ or $\Omega_h^n$ are assigned a nodal value using the Lagrangian interpolation (e.g., marked by \textcolor{gray}{\tiny $\blacksquare$} in the top row of Figure~\ref{fig:nodal_values_of_common_triangulation}).
	\end{enumerate}
	
\paragraph{Interpolation errors of common interpolation operators.}
For later analysis it will be essential to quantify the interpolations errors of the common interpolation operators $\mathcal{I}_k^{\comesh}$. 
We write $E_{\Ga_h}^{\comesh}$ as the set of boundary edges of $\mathcal{T}_h^{\comesh}$.
Additionally, we split the element set disjointly into the set of interior elements $\mathcal{T}_{h,\text{int}}$, with at most one vertex on the boundary, and the already introduced set of boundary elements $B_h$, see Lemma~\ref{Lemma:GeomAppError}, such that $\mathcal{T}_{h,\text{int}} \dot{\bigcup} B_h = \mathcal{T}_h$.

Then the following holds
\begin{lem}
	\label{Lemma:common_interpolation}
	Let $u_h^n = (\phi_h^{n},\psi_h^{n}) \in \calV_{h}^n$ and let the common domains be defined as above, in particular we assume refinement and coarsening by NVB with immediate lifting, see Remark~\ref{remark:NVB}. 
	For each element $T \in \mathcal{T}_{h,\text{int}}^{n}$, or for each element $T \in B_h^{n}$ whose boundary edge $S$ fulfills $S \in E_{\Ga_h}^{n} \cap E_{\Ga_h}^{\comesh}$, i.e. the element is not effected by lifting associated to refinement, we have
	\begin{align}
	\label{eq:common_int_zero}
		\|(u_h^n)^{\ell^n}-(\mathcal{I}_n^{\comesh} u_h^n)^{\ell^\comesh}\|_{H^1(T^\ell)} = 0  
	\end{align}
	If on the other hand $T \in B_h^n$ and the boundary edge fulfills $S \notin E_{\Ga_h}^{n} \cap E_{\Ga_h}^{\comesh}$, i.e. the element was atleast refined once on a boundary element, and thus a vertex was lifted, to derive the child elements in the common triangulation, then
	\begin{equation}
	\label{eq:common_int_bulk}
		\begin{aligned}
		\|(\phi_h^n)^{\ell^n}- (\mathcal{I}_n^{\comesh} \phi_h^n)^{\ell^\comesh}\|_{L^2(T^\ell)} &\leq c h_{T}^2 |\phi_h^n|_{H^1(T)}, \\
		|(\phi_h^n)^{\ell^n}- (\mathcal{I}_n^{\comesh} \phi_h^n)^{\ell^\comesh}|_{H^1(T^\ell)} &\leq c  h_{T} |\phi_h^n|_{H^1(T)},
		\end{aligned}
\end{equation}
	and on the corresponding boundary edge $S = T \cap \Gamma_h^n$ on the surface we have the following interpolation error
	\begin{equation}
	\label{eq:common_int_surf}
		\begin{aligned}
			\|(\psi_h^n)^{\ell^n}- (\mathcal{I}_n^{\comesh} \psi_h^n)^{\ell^\comesh}\|_{L^2(S^\ell)} &\leq c  h_{\tilde{S}}^2 |\psi_h^n|_{H^1(S)}, \\
			|(\psi_h^n)^{\ell^n}- (\mathcal{I}_n^{\comesh} \psi_h^n)^{\ell^\comesh}|_{H^1(S^\ell)} &\leq c  h_{\tilde{S}} |\psi_h^n|_{H^1(S)}.
		\end{aligned} 
	\end{equation} 
	The scalar $c$ is an $h-$ and $n-$uniform constant, depending only on $\Gamma$ and $d$. Further the same bounds hold for $(\psi_h^{n-1},\psi_h^{n-1}) \in \calV_h^{n-1}$ with the interpolant $\mathcal{I}_{n-1}^{\comesh}$ and the corresponding parent meshes $\Om_h^{n-1}$ and $\Ga_h^{n-1}$.
\end{lem}

\begin{proof}
First note that a proof for the results for boundary edges \eqref{eq:common_int_surf} is given in \cite{KovacsLantelme2026}.
In comparison to \cite{KovacsLantelme2026} we avoid the dependence of the constant on the refinement level by representing the interpolation errors on the parent meshes as done in \cite{Lantelme2026}.
All upcoming argument solely focus on the analysis of the interpolation error in the bulk, where we follow the proof concepts for coarsening interpolation errors on surfaces, described in \cite{KovacsLantelme2026}, and aim to extend the results to bulk domains.

\textbf{Interior elements and results for surface domains}
By construction of $\Omega_{h}^{\comesh}$ the nodal values of vertices present in both $\mathcal{T}_{h}^{\comesh}$ and $\mathcal{T}_{h}^{n}$ are identical. Further in the bulk, even if $\Omega_{h}^{n}$ is a refinement of $\Omega_{h}^{\comesh}$ the finite element function of $(\phi_h^{n},\psi_h^{n}) \in \calV_h^{n}$ exactly matches the finite element function of $(\mathcal{I}_n^{\comesh} \phi_h^{n},\mathcal{I}_n^{\comesh} \psi_h^{n}) \in \calV_h^{\comesh}$ for all unchanged elements, in particular for all interior elements $T \in \mathcal{T}_{h,\text{int}}^{n}$. Note that this only works as we work with linear finite elements and our interpolation is invariant on elements which are unaffected by lifting. Similarly, the different lifts align for unchanged elements, thus \eqref{eq:common_int_zero} holds.
Additionally, boundary elements which are not refined by a bisection of a boundary edge are not effected by intermediate node-lifting and thus do not introduce geometric interpolation errors. 
 
Thus the focus of our analysis are the boundary elements where at least one refinement occurred on a boundary edge.
To tackle the problem we assume that our element $T$ was only refined, with a new vertex on the boundary, utilizing newest vertex bisection, but generalization to multiple refinements yields identical bounds \cite{KovacsLantelme2026}. Using the newest vertex bisection guarantees that at most one vertex is lifted in each element, we will index the $l$-th vertex, i.e. the final boundary vertex as the lifted vertex.  
In comparison to the surface case the main difference is that we have to extend the results to account for the bulk lift and the underlying bulk domain.

\textbf{Transformations via $\theta$-arguments}
Using the common triangulation $\Omega_h^{\comesh}$ (see Figure~\ref{fig:nodal_values_of_common_triangulation}) as a basis, we use a $\theta$-argument very similar to the one developed in Section~4 of \cite{KLLP_2017}. 
	
	To this end we introduce a continuous affine map which transforms a function from the discrete mesh $\Omega_h^0 := \widehat{\Omega}_h^{n}$ (with an unlifted refined node, see Figure~\ref{fig:nodal_values_of_common_triangulation} top row), to its common interpolated version  $\Omega_h^1 := \Omega_h^{\comesh}$ with lifted nodes by construction of the common triangulation (see Figure~\ref{fig:nodal_values_of_common_triangulation} bottom row).
	We define the $\theta$-dependent discrete surface $\Omega_h^\theta$ elementwise by the refined boundary elements for $\theta \in [0,1]$, such that each moving node $x_l =: x_l^0$ is moved linearly to $x_l^1 := (x_l)^\ell$ via the unique point projection $p$, i.e., $x_l^\theta := x_l^0 + \theta (x_l^1 - x_l^0)$.  Based on the $\theta$-dependent mesh we define finite element spaces $\calV_h^\theta$.
	
	Similar to Section~\ref{sec:lift_interpolation} we define the following quantities pointwise in $\theta$ for the bulk domain $\Omega_h^\theta$ 
	We define the following objects:
	\begin{center}
		\begin{tabular}{  m{3cm}  m{9.2cm} } 
			$T^\theta$&  element of $\mathcal{T}_h^\theta$\\  
			$F_T^\theta \colon \hat T \times \theta \to T^\theta$ & affine reference map $F_T^\theta(\hat x) = A(\theta) \hat x + b$ for $\hat x \in \hat T$ \\
			$\Phi^\theta \colon \hat T \times \theta \to \R^d $ & bulk lift component $\Phi^\theta(\hat{x}) = (\hat{\lambda}^\star)^{3} (p(y^\theta)-y^\theta)$ \\			
			$F_T^{e,\theta} \colon \hat T \to (T^\theta)^{\ell^\theta}$ & reference map of a curved simplex, given by $F_T^{e,\theta} = F_T^\theta + \Phi^\theta$ \\
			$G_h^\theta \colon T^\theta \to (T^\theta)^\ell$ &  lift map of the element $T$ given by $F_T^{e,\theta} \circ (F_T^\theta)^{-1}$ \\
			$(\cdot)^{\ell^\theta} \colon \Omega_h^\theta \rightarrow \Omega$ &	lift operator \\	
			\makecell[l]{$x_1,\dots,x_{l-1},x_l^\theta,$\\$x_{l+1},\dots,x_{d+1}$} &  nodes of an element $T^\theta$, where the, potential empty, set of boundary nodes are always labeled first, i.e. $x_i \in \Ga_h^n$ for $i = 1,\dots,l$, where $l$ denotes the number of boundary nodes.
		\end{tabular}
	\end{center}
	We want to highlight the key dependencies on $\theta$. First note that we assume that only a single vertex moves, which is labeled $x_l^\theta$, all other are fixed.
	This also justifies why the affine reference map $F_T^\theta$ can be written with a fixed $b$, which we will later choose as $b = x_1$.
	For the bulk-lift component $\Phi^\theta$ we note that the barycentric coordinates $\hat{\lambda}_i$ are independent of $\theta$ but $y^\theta$ depends on the node positions, in particular on $x_l^\theta$.
	
	If we represent an arbitrary point $x^0 \in T^{0} \subset \Omega_h^0$ by its barycentric coordinates $\lambda_i$, the pushed-forward point $x^\theta \in T^\theta$ can be written by $x^\theta = \lambda_l x_l^\theta + \sum_{i\neq l} \lambda_i x_i $ by construction. Note that $\lambda_i$ are independent of $\theta$.

	Now a function in $V_h^\theta$, the finite element space corresponding to $\Omega_h^\theta$, with fixed nodal vector $(w_i)$ can be written as $w_h^\theta = \sum_i w_i \varphi_i^\theta$ with fixed nodal values $w_i$ and the nodal basis functions $\varphi_i ^\theta \in V_h^\theta$.
	An important property is the following relation between $\varphi_i ^\theta$ and the reference basis functions $\hat{\varphi}_{\hat{i}}$, with $\hat{i}$ the local vertex numbering on the reference simplex $\hat{T}$.
	Note for $x \in T^\theta$ we have the relation $\varphi_i ^\theta(x) = \hat{\varphi}_{\hat{i}}(\hat{x}) := \hat{\varphi}_{\hat{i}}((F^\theta)^{-1} (x))$.

\textbf{General approach by fundamental theorem of calculus}
The key idea is now to rewrite the difference in \eqref{eq:common_int_bulk} by the fundamental theorem of calculus and careful differentiation.
The upcoming arguments are always understood elementwise, in particular the maps related to the lift and affine transformations are defined on the related element of some point. 
We introduce the notation $x^\theta(a) := (G_h^\theta)^{-1}(a)$ for the unlifted point of $a \in \Omega$ onto $\Omega_h^\theta$.

For any point $a \in \Omega$ we rewrite the difference of two lifted functions as follows (cf. Lemma~4.1 in \cite{KLLP_2017}):
	\begin{equation}
	\label{eq:coarsening_proof_thm_of_calculus}
		\begin{aligned}
			(\phi_h^1)^{\ell^1}(a)-(\phi_h^0)^{\ell^0}(a) 
			= &\ \int_0^1 \frac{\d}{\d \theta}(\phi_h^\theta)^{\ell^\theta}(a) \d \theta 
			= \int_0^1 \frac{\d}{\d \theta} \Big(\phi_h^\theta \big(x^\theta(a)\big) \Big) \d \theta \\
			= &\ \int_0^1  \big(D_x \phi_h^\theta\big)\big(x^\theta(a)\big) \cdot \frac{\d}{\d \theta} x^\theta(a) + (\partial_\theta \phi_h^\theta)\big(x^\theta(a)\big)  \d \theta  \\
			= &\ \int_0^1\big(D_x \phi_h^\theta\big)\big(x^\theta(a)\big) \cdot \Big(\frac{\d}{\d \theta} x^\theta(a) - \big(\partial_\theta F^\theta\big)\big((F^\theta)^{-1}\big(x^\theta(a)\big) \Big) \d \theta .
		\end{aligned}
	\end{equation}

	Where the last line follows from noting that by the transport property  
\begin{align*}
	0 = \frac{\d}{\d\theta} \phi_h^\theta (F^\theta(\hat{x}))= (D_x \phi_h^\theta)(F^\theta(\hat{x})) \cdot \partial_\theta F^\theta(\hat{x}) + (\partial_\theta \phi_h^\theta)(F^\theta(\hat{x})),
\end{align*}	
	and evaluating in $\hat{x} = (F^\theta)^{-1}(x^\theta(a))$.

	In some sense the difference in \eqref{eq:coarsening_proof_thm_of_calculus} measures the $\theta$-dependent motion of the point $x^\theta(a)$, first due the explicit dependence in the lift map $G^\theta$ and then for the reference map $F^\theta$. 

Similar we can write for the gradients,
\begin{equation}
	\begin{aligned}
		\label{eq:coarsening_proof_stiffness_fdm_theorem}
		\nabla (\phi_h^1)^{\ell^1}(a)&-\nabla (\phi_h^0)^{\ell^0}(a) =  D_a \int_0^1\frac{\d}{\d \theta} \Big(\phi_h^\theta (x^\theta(a)) \Big) \d \theta\\
		 = &\int_0^1 \big(D_x \phi_h^\theta\big)\big(x^\theta(a)\big) \cdot D_a \Big(\frac{\d}{\d \theta} x^\theta(a)- \big(\partial_\theta F^\theta\big)\big((F^\theta)^{-1}(x^\theta(a)) \Big) \d \theta.
		 \end{aligned}
	\end{equation}
Note that after swapping differentiation and integration we used that one part of the product rule vanishes elementwise due to the linearity of the finite element function $\phi_h^\theta$, since
	\begin{align*}
	D_a \big[\big(D_x \phi_h^\theta\big)\big(x^\theta(a)\big)\big] = \big(D_x^2 \phi_h^\theta\big) \big(x^\theta(a)\big) D_a x^\theta(a) = 0.
	\end{align*}

\textbf{Bounding the derivatives}
To show the Lemma we bound for each parent element $T \in \Omega_h^n$:
\begin{align}
\Big\|\frac{\d}{\d \theta} x^\theta - \big(\partial_\theta F^\theta\big)\big((F^\theta)^{-1}\big(x^\theta\big)\Big\|_{L^\infty(T^{\ell})} &\leq ch_{T}^2 \\
\Big\|D_a \Big(\frac{\d}{\d \theta} x^\theta- \big(\partial_\theta F^\theta\big)\big((F^\theta)^{-1}(x^\theta \Big)\Big\|_{L^\infty(T^{\ell})} &\leq ch_{T}.
\end{align}
These estimates are proven by careful computations in Appendix~\ref{sec:app_estimates_int_bounds}.

We can identify the gradient as a lifted quantity and based on norm equivalence type results in the bulk by Remark~\ref{rem:normequivalence} we can show
	\begin{equation}
		\label{eq:coarsening_proof_norm_eq}
		\|\big(D_x \phi_h^\theta\big)\big(x^\theta\big)\|_{L^2((T^\theta)^\ell)} = \|(D_x \phi_h^\theta)^{\ell^\theta} \|_{L^2((T^\theta)^\ell)} \leq c |\phi_h^\theta|_{H^1(T^\theta)}.  
	\end{equation}	

The $\theta$-uniformity of $\phi_h^\theta$ \cite{KLLP_2017} then further implies 
$$|\phi_h^\theta|_{H^1(\Omega_h^\theta)} \leq c |\phi_h^0|_{H^1(\Omega_h^0)} = c |\phi_h^n|_{H^1(\Omega_h^n)},$$
which yields the interpolation error bounds of \eqref{eq:common_int_bulk}.

The main idea to justify multiple refinement/coarsening is to introduce intermediate meshes given by refining/coarsening one level at a time. The arguments are applied for consecutive meshes. Afterwards by $\theta$-uniformity and carefully tracking the element-size, we are able to combine a bound on the parent element, as done in \cite{Lantelme2026}. 
\end{proof}



\section{Reliability analysis for $K,L \in (0,\infty)$} \label{sec:Reliability}

This section is devoted to the reliability analysis of our adaptive method by an a posteriori error estimate for the error indicators. This is achieved by first relating the error and the residual, and then bounding the residual by appropriate error indicators.

\subsection{Reduction to quadratic nonlinear potentials}

The nonlinear potentials $F_\Om$ and $F_\Ga$ are of the form: $F\in C^2(\R)$ such that $F(r)=F_1(r) - \frac\lambda2 r^2$, $F$ is bounded from below, $F_1$ is non-negative and convex and has polynomial growth, i.e. there exists $q \in [1,\infty)$ such that, for some $a>0$, 
\begin{align}
	|F_1'(r)|\leq a|r|^q+a.\label{eq:RegularityF}
\end{align}
An important example of this form is the double-well potential $F(r)=\frac{(1-r^2)^2}4 $.

In the following we only consider quadratic potentials $F_\Om$ and $F_\Ga$, i.e. $q=1$, and hence with a Lipschitz-continuous derivative $F_\Om'$ and $F_\Ga'$. 

This case is generalized to the proposed regularity \eqref{eq:RegularityF} by an analogous approach as in \cite[Section~4.2]{ElliottSales2024}: We argue that in the general case $q>1$ we have by the regularity of the exact solution $u=(\phi,\psi)$:
\begin{align*}
	\max_{t\in[0,T]}\max_{x \in \Om}|\phi(x,t)| + \max_{t\in[0,T]}\max_{x \in \Ga}|\psi(x,t)| \leq u_{\max}.
\end{align*}
Therefore we instead work with the quadratic potentials $F_{\Om,M} \in C^2(\R)$, defined by
\begin{alignat}{2}
	\label{eq:quadratic_potential}
	F_{\Om,M}(r):=\begin{cases}
		F_\Om(M) + F_\Om'(M)(r-M)+\frac{F_\Om''(M)}{2} (r-M)^2, \quad &r>M \\
		F_\Om(r), \quad &|r|\leq M \\
		F_\Om(-M) + F_\Om'(-M)(r+M)+\frac{F_\Om''(-M)}{2} (r+M)^2, \quad &r<-M
	\end{cases},
\end{alignat}
and analogously $F_{\Ga,M} \in C^2(\R)$, for some $M>u_{\max}$. What remains to be seen is that we can choose $M$, so that $$ |\phi_{h}^n|+|\psi_{h}^n| \leq M, $$ independent of $h$ and $\tau$. This can be established by an optimal-order convergence result similar to \cite[Theorem~5.1]{Bullerjahn2025} and an inverse estimate for finite element norms, and comes with a step-size restriction $\max_{n} \tau_n \leq Ch^{1+\eps}$ for some arbitrarily chosen positive constants $C,\eps$.

\subsection{Bounding the error by the residual}
We consider the errors $e_u(t)=u(t)-u_{h,\tau}(t)$ and $e_w(t)=w(t)-w_{h,\tau}(t)$ of the exact solution and the numerical approximation, see \eqref{eq:u_cont} and \eqref{eq:Def w_htau}. Then the errors satisfy the error equations
\begin{subequations} \label{eq:ErrorEq}
	\begin{align}
		m(\partial_t e_u,\zeta) + a^{L,\beta}(e_w,\zeta) =&\ -m(R^1,\zeta), \label{eq:ErrorEq1}\\
		m(e_w,\eta) - a^{K,\alpha,\eps}(e_u,\eta) =&\  m(\inveps r,\eta) - m(R^2,\eta),\label{eq:ErrorEq2}
	\end{align}
\end{subequations}
for all $\zeta \in \calV$, $\eta \in \calV$, and where $r=F'(u)-F'(u_{h,\tau})$. 

Here the residual terms are defined as
\begin{subequations}
	\begin{align}
		m(R^1,\zeta)=&\ m(\partial_t u_{h,\tau},\zeta) + a^{L,\beta}(w_{h,\tau},\zeta), \label{eq:DefResidual1} \\
		m(R^2,\eta)=&\ m(w_{h,\tau},\eta) - a^{K,\alpha,\eps}(u_{h,\tau},\eta) - m(\inveps F'(u_{h,\tau}),\eta), \label{eq:DefResidual2}
	\end{align}
\end{subequations}
for all $\zeta \in \calV$ and $\eta \in \calV$.

An important tool in the following analysis is the problem specific inverse Laplace operator: We define the operator $\big(-\Delta_{L,\beta}\big)^{-1} \colon \calH \to \calV_0$ for $\phi \in \calH$, with the solution space $\calV_0 := \{v \in \calV \colon \langle v \rangle =0\}$, by the solution to the elliptic problem
\begin{align}
	a^{L,\beta}(\big(-\Delta_{L,\beta}\big)^{-1} (\phi),\zeta) = m(\phi - \langle \phi \rangle (\beta,1),\zeta), \forallquad \zeta \in \calV, \label{eq:Def invLaplace}
\end{align}
where we denote for $\phi = (\zeta,\eta) \in \calH$ its mean by 
$$\langle \phi \rangle := \frac{\beta \int_\Om \zeta + \int_\Ga\eta}{\beta^2 |\Om| + |\Ga|} .$$
Then an application of the Lax-Milgram theorem shows that this is well defined and satisfies
\begin{align}
	\| \big(-\Delta_{L,\beta}\big)^{-1} \phi \|_\calV \leq C \|\phi \|_\calH, \label{eq:Stability inv laplace}
\end{align}
similar to \cite[Proposition~A.1]{KnopfStange2024}. Using this operator we define the $\calH_{L,\beta}^{-1}$-seminorm as
\begin{align}
	| \phi |_{\calH_{L,\beta}^{-1}}^2 := m(\big(-\Delta_{L,\beta}\big)^{-1} \phi,\phi), \forallquad \phi \in \calV. \label{eq:def H-1 norm}
\end{align}

\begin{lem}
	\label{Lemma:error_residual_bound}
	Let $(u,w)$ be a solution of the bulk--surface Cahn--Hilliard equation with dynamic boundary conditions \eqref{eq:CHweak}, with quadratic potentials, i.e. \eqref{eq:RegularityF} with $q=1$. Then the following bound holds, connecting the error and the residual:
	\begin{align}
		|e_u|_{L^\infty(0,t^n;\calH_{L,\beta}^{-1})}^2 + \eps \|e_u\|_{L^2(0,t^n;\calV)}^2 \leq&\ C \exp\big(c\eps^{-3} t^n\big) \Big( R_h^n+ \calE_I \Big),\label{eq:ErrorByResidualBound}
	\end{align}    
	where $c$ and $C$ are independent of $h$, $\tau$ and $\eps$, but depends on the Lipschitz-constant of $F'$, and the residuals are collected in 
	\begin{align}
		R_h^n:= \eps \| R^1\|_{L^2(0,t^n;\calV^{\ast})}^2 + \inveps \| m(R^1,(\beta,1))\|_{L^2(0,t^n)} + \eps^{-1} \|R^2\|_{L^2(0,t^n;\calV^{\ast})}^2 , \label{eq:DefResidualNorms}
	\end{align}
	and the error in the initial data in
	\begin{align}
		\calE_I := | u^0-(u_h^0)^\ell |_{\calH_{L,\beta}^{-1}}^2 + \inveps |m(u^0-(u_h^0)^\ell,(\beta,1))|^2. \label{eq: Def Initial error}
	\end{align}
\end{lem}

\begin{proof}
	In the following $c$ is a positive generic constant, independent of $h$ and $\tau$ that takes different values on different occurrences. By $\rho >0$ we denote a small number, independent of $h$ and $\tau$, used in Young's inequality, and hence we will often incorporate multiplicative constants into those yet unchosen factors.
	
	We start by testing the first equation \eqref{eq:CHweakD1} with $\big(-\Delta_{L,\beta}\big)^{-1} e_u$ and subtract the second equation \eqref{eq:CHweakD2} tested with $e_u$, to obtain
	\begin{align}
		\ddt \frac12 |e_u|_{\calH_{L,\beta}^{-1}}^2 + | e_u |_{a^{K,\alpha,\eps}}^2 =&\ -m(R^1,\big(-\Delta_{L,\beta}\big)^{-1} e_u) - m(\inveps r, e_u) \nonumber \\
		&\ + m(R^2,e_u)\nonumber \\
		=:&\ (\text{RHS}).\label{eq: energy tested}
	\end{align}
	
	Then we estimate the right-hand side by \eqref{eq:Stability inv laplace}, Cauchy-Schwarz and Young's inequality, the Lipschitz continuity of $F'$, and the definition of the dual norm, to
	\begin{align}
		(\text{RHS}) \leq&\ \delta \eps \| e_u \|_{\calV}^2 + c \inveps \| e_u \|_\calH^2 + c \eps \| R^1 \|_{\calV^\ast}^2 + c \inveps \| R^2 \|_{\calV^\ast}^2. \label{eq:energy RHS 1}
	\end{align}
	Then we obtain by \eqref{eq:Def invLaplace},
	\begin{align}
		\inveps \| e_u \|_\calH^2 =&\ \inveps m(e_u,e_u) = \inveps a^{L,\beta}\big((-\Delta_{L,\beta})^{-1} e_u, e_u \big) + \inveps \frac{|m( e_u,(\beta,1)|^2}{ \beta^2 |\Om| + |\Ga|} \nonumber \\
		\leq&\ \inveps| (-\Delta_{L,\beta})^{-1} e_u|_{a^{L,\beta}} |e_u |_{a^{L,\beta}} + c \inveps|m(e_u,(\beta,1))|^2 \nonumber \\
		\leq&\ c \inveps |e_u|_{\calH_{L,\beta}^{-1}} \|e_u \|_{\calV} + c \inveps|m(e_u,(\beta,1))|^2 \nonumber \\
		\leq&\ \delta \eps \| e_u \|_\calV^2 + c \eps^{-3} |e_u|_{\calH_{L,\beta}^{-1}}^2 + c \inveps|m(e_u,(\beta,1))|^2. \label{eq:energy RHS 2}
	\end{align}

	Inserting now the estimates \eqref{eq:energy RHS 1} and \eqref{eq:energy RHS 2} into \eqref{eq: energy tested}, we obtain
	\begin{align}
	\label{eq:energy before poin}
		\ddt \frac12 |e_u|_{\calH_{L,\beta}^{-1}}^2 + | e_u |_{a^{K,\alpha,\eps}}^2 \leq&\ \delta \eps \| e_u \|_\calV^2  + c \eps^{-3} |e_u|_{\calH_{L,\beta}^{-1}}^2 + c \inveps|m(e_u,(\beta,1))|^2 \nonumber \\
		&\ + c \eps \| R^1 \|_{\calV^\ast}^2 + c \inveps \| R^2 \|_{\calV^\ast}^2.
	\end{align}
	
	Next we derive an \emph{almost mass conservation} of the error by testing the first equation with the constant function $(\beta,1) \in \calV$ and integrate to obtain
	\begin{align}
		|m(e_u(t),(\beta,1))|^2 \leq&\ c |m(e_u(0),(\beta,1))|^2 + c \int_{0}^{t^n} |m(R^1,(\beta,1))|^2 . \label{eq:Almost mass conservation}
	\end{align}
	This leads to a Poincaré--Wirtinger-type estimate using the bulk--surface Poincaré-inequality, see e.g. \cite[Lemma~A.1]{KnopfLiu2021}:
	\begin{align*}
		\|e_u(t)\|_{V}^2 \leq c |e_u(t)|_{a^{K,\alpha,1}}^2 + c |m(e_u(0),(\beta,1))|^2 + c \int_{0}^{t^n} |m(R^1,(\beta,1))|^2. 
	\end{align*}
	Multiplying with $\eps<1$, we finally obtain
	\begin{align}
		\eps \|e_u(t)\|_{V}^2 \leq c |e_u(t)|_{a^{K,\alpha,\eps}}^2 + c \eps |m(e_u(0),(\beta,1))|^2 + c \eps \int_{0}^{t^n} |m(R^1,(\beta,1))|^2. \label{eq:PoincareU}
	\end{align}
	
	Integrating \eqref{eq:energy before poin}, inserting \eqref{eq:Almost mass conservation} and \eqref{eq:PoincareU}, and absorbing terms we arrive at
	\begin{align*}
		|e_u(t^n)|_{\calH_{L,\beta}^{-1}}^2 + \eps \|e_u\|_{L^2(0,t^n;\calV)}^2 \leq&\ c \eps^{-3} | e_u |_{L^2(0,t^n;\calH_{L,\beta}^{-1})}^2 + c \Big( R_h^n + \calE_I \Big),
	\end{align*}
	and an application of Gr\"onwall's inequality yields the final result.
\end{proof}

\subsection{Notational conventions}

In the following section we deal with the estimation of the residual term $R_h^n$ by error estimators $\eta^n$. So far we have been sticking to a notation for the combined bulk--surface, but as it turns out these estimates force us to distinguish between the bulk and surface. Nevertheless, we want to stick to a bulk--surface notation as this lets us see the overall structure of the estimates in a clearer way, and reduces the amount of terms. 

An example for a naturally appearing term and its bulk surface notation, for $v_h^n=(\zeta_h^n,\eta_h^n) \in \calV_h^n$, is the following:
\begin{equation*}
	\sum_{(T,S) \in \calT_h^n \times E_{\Ga_h}^n} (h_T^{k_1},h_S^{k_2}) \| v_h^n \|_{m_h^{n,T,S}}^2 := \sum_{T \in \calT_h^n} h_T^{k_1} \| \zeta_h^n \|_{L^2(T)}^2 + \sum_{S \in E_{\Ga_h}^n } h_S^{k_2} \| \eta_h^n \|_{L^2(S)}^2,
\end{equation*}
which includes the notation by the inner product $m_h^{n,T,S}$ for the $L^2(T) \times L^2(S)$-norm. Whenever possible, we further simplify $(h_T^k,h_S^k)$ to $h^k$. Analogously, we define this, for example, for the semi-inner product $a_h^{K,\alpha,\eps}$, and in some cases with $\nabla v_h^n = (\nabla \zeta_h^n,\nabla_{\Ga_h^n} \eta_h^n)$ or $[\nabla v_h^n \cdot \nu]:=([\nabla \zeta_h^n \cdot \nu_S],[\nabla_{\Ga_h^n}\eta_h^n \cdot \nu_P])$. We also shorten the notation $\comesh$ to $\oplus$, as all our considerations are on a single time interval, hence the common mesh does not change.

Using this notation we are able to stick to a bulk--surface notation for almost all terms, except in some instances for $a_h^{K,\alpha,\eps}$, because of the substantially different estimates for the stiffness terms and the boundary coupling terms, c.f. \eqref{eq:Def bilinear forms}. In these cases we always denote the numerical solution as $u_h^n=(\phi_h^n,\psi_h^n)$ and $w_h^n=(\mu_h^n,\theta_h^n)$.

\subsection{Estimating the residuals $R^1$ and $R^2$}
\label{sec:indicators_bound_residual_bound}
In order to estimate the first term in \eqref{eq:DefResidualNorms}, we use standard techniques for residual based a posteriori error estimates, see e.g. \cite{Verfuerth2003} and \cite[Chapter~9]{BrennerScott2008}, as well as geometric error bounds, see \cite{ElliottRanner2013}, and Lemma~\ref{Lemma:common_interpolation}. The resulting error indicators for $n=1,\dots,N$ are then given by
\begin{subequations} \label{eq:Error indicators 1}
	\begin{align}
		\eta_1^n :=&\ \eta_{h,1}^n + \eta_{c,1}^n + \eta_{M,1}^n,\\
		(\eta_{h,1}^n)^2 :=&\ \sum_{(T,S) \in \calT_h^n \times E_{\Ga_h}^n} h^2 \| \partial_{\tau_n} u_h^n \|_{m_h^{n,T,S}}^2 + \sum_{(S,P) \in E_{\text{int}}^n \times P_{\Ga_h}^n} h \| [\nabla w_h^n \cdot \nu]\|_{m_h^{n,S,P}}^2 \nonumber \\
		&\ + \sum_{S \in E_{\Ga_h}^n} h_S \| \partial_\nu \mu_h^n - L^{-1} (\beta \theta_h^n - \mu_h^n)\|_{L^2(S)}^2 + h_S^2 \| L^{-1}(\beta \theta_h^n - \mu_h^n) \|_{L^2(S)}^2, \\
		(\eta_{c,1}^n)^2 :=&\  \sum_{(T,S) \in \calT_h^\oplus \times E_{\Ga_h}^\oplus} \tau_n^{-2} \| \calI_{n-1}^\oplus u_h^{n-1} - \calI_{n}^\oplus \calI_{n-1}^n u_h^{n-1}\|_{m_h^{\oplus,T,S}}^2 ,\\
		(\eta_{M,1}^n)^2 :=&\  \sum_{(T,S) \in \calT_h^\oplus \times E_{\Ga_h}^\oplus} \tau_n^{-2} | m_h^{\oplus,T,S}\big(\calI_{n-1}^\oplus u_h^{n-1} - \calI_{n}^\oplus \calI_{n-1}^n u_h^{n-1},(\beta,1)\big)|^2,
	\end{align}
	with $E_{\text{int}}^n := E_h^n \setminus E_{\Ga_h}^n$ and the geometric estimators given by
	\begin{align}
		\calG_1^n :=&\ \calG_{g,1}^n + \calG_{c,1}^n, \\
		(\calG_{g,1}^n)^2 :=&\ \sum_{(T,S) \in B_h^n \times E_{\Ga_h}^n} (h_T^2,h_S^4) \| \nabla w_h^n \|_{m_h^{n,T,S}}^2 + (h_T^2,h_S^4) \| \partial_{\tau_n} u_h^n \|_{m_h^{n,T,S}}^2 \nonumber \\
		&\ + \sum_{(T,S) \in B_h^\oplus \times E_{\Ga_h}^\oplus} (\frac{h_T^2}{\tau_n^2},\frac{h_S^4}{\tau_n^2}) \| \calI_{n-1}^\oplus u_h^{n-1} - \calI_{n}^\oplus \calI_{n-1}^n u_h^{n-1} \|_{m_h^{\oplus,T,S}}^2, \\
		(\calG_{c,1}^n)^2:=&\ \sum_{(T,S) \in \hat B_h^{n-1} \times \hat E_{\Ga_h}^{n-1}} \frac{h^4}{\tau_n^{2}}  \| \nabla u_h^{n-1} \|_{m_h^{n-1,T,S}}^2 \nonumber \\
		&\ + \sum_{(T,S) \in \hat B_h^{n} \times \hat E_{\Ga_h}^{n}} \frac{h^4}{\tau_n^{2}}  \| \nabla  (\calI_{n-1}^n u_h^{n-1}) \|_{m_h^{n,T,S}}^2 
	\end{align}
\end{subequations}
where we recall that $B_h$ collects elements which have at least two nodes on the boundary and the element set $\hat B$ and edge set $\hat E_{\Ga_h}$ refer to the subset of elements, respectively boundary edges, which have at least one boundary vertex which was either refined or coarsened when solving for $n$.

\begin{lem} \label{Lemma:Residual bound 1}
	Let $n=1,\dots,N$ and $t \in (t^{n-1},t^n]$, then the dual norm \eqref{eq:dualNormDef} of the first residual $R^1$ fulfills
	\begin{align*}
		\| R^1(t) \|_{\calV^{\ast}} \leq&\ C \Big(\eta_1^n + \calG_1^n \Big), \\
		|m(R^1(t),(\beta,1))|  \leq&\ C \Big( \eta_{M,1}^n + \calG_{1}^n\Big)
	\end{align*}
	for some constant $C>0$ independent of $h$, $\tau$ and $\eps$, and the error indicators defined in \eqref{eq:Error indicators 1}.
\end{lem}

\begin{proof}
	We start with the dual norm and by decomposing the residual by successively adding zeros, to obtain with \eqref{eq:u_cont}, \eqref{eq:u_computable} and \eqref{eq:Def w_htau}:
	\begin{equation} \label{eq:Res 1 split}
		\begin{aligned}
			m(R^1,\zeta) =&\ m(\partial_t u_{h,\tau}, \zeta) + a^{L,\beta}(w_{h,\tau},\zeta) \\
			=&\ m(R_c^1,\zeta) + m(R_g^1,\zeta) + m(R_h^1,\zeta), 
		\end{aligned}
	\end{equation}
	with the coarsening, geometric and spatial residuals
	\begin{align*}
		m(R_c^1,\zeta):=&\ \frac1{\tau_n} m((\calI_{n-1}^n u_h^{n-1})^{\ell_n} - (u_h^{n-1})^{\ell_{n-1}}, \zeta), \\
		m(R_g^1,\zeta):=&\ \big[m((\partial_{\tau_n} u_h^n)^{\ell_n}, \zeta) - m_h^n(\partial_{\tau_n} u_h^n, \zeta^{-\ell})\big] \nonumber \\
		&\ + \big[a^{L,\beta}((w_h^n)^{\ell_n},\zeta) - a_h^{n,L,\beta}(w_h^n,\zeta^{-\ell})\big], \\
		m(R_h^1,\zeta):=&\ m_h^n(\partial_{\tau_n} u_h^n, \zeta^{-\ell}) + a_h^{n,L,\beta}(w_h^n,\zeta^{-\ell}).
	\end{align*}
	Here, we have used the fact that the linear interpolation exactly satisfies $\partial_t \widetilde u_{h,\tau} = \partial_{\tau_n} u_h^n$.
	
	The spatial residual $R_h^1$ resembles exactly the fully discrete equation \eqref{eq:CHweakD1}, the only difference being that we now test with a general $\zeta^{-\ell_n}=(\eta^{-\ell_n},\xi^{-\ell_n}) \in H^1(\Om_h^n) \times H^1(\Ga_h^n)$. As the discrete equation holds for all functions from our finite element space $\calV_h^n$ we can insert a finite element function given by the tuple $\zeta_h = (\eta_h,\xi_h) \in \calV_h^n$ of Scott-Zhang interpolant on the bulk and on the surface of $\zeta^{-\ell_n}$, see, e.g. \cite{BrennerScott2008,CamachoDemlow2015} for their definitions and interpolation estimates. The estimate by the spatial estimator is then obtained by splitting the bilinear forms into bulk and boundary parts:
	\begin{align}
		m(R_h^1(t),\zeta) =&\ \int_{\Om_h^n} \partial_{\tau_n} \phi_h^n (\eta^{\ell_n} - \eta_h) + \int_{\Ga_h^n} \partial_{\tau_n} \psi_h^n (\xi^{\ell_n} - \xi_h) \nonumber \\
		&\ + \int_{\Om_h^n} \nabla \mu_h^n \nabla(\eta^{\ell_n} - \eta_h) + \int_{\Ga_h^n} \nabla_{\Ga_h^n}\theta_h^n \nabla_{\Ga_h^n}(\xi^{\ell_n} - \xi_h) \nonumber \\
		&\ + L^{-1} \int_{\Ga_h^n} (\beta\theta_h^n-\mu_h^n) \big(\beta(\xi^{\ell_n} - \xi_h) - (\eta^{\ell_n} - \eta_h)\big). \label{eq:R1 spatial splitted}
	\end{align}
	We follow the typical arguments (see, e.g. \cite[Chapter~9]{BrennerScott2008}): Divide the domain into its elements; the stiffness term is first integrated by parts, where we can simplify the terms by introducing the jump. For the stiffness term in the bulk $\Om_h^n$ the introduction of the jump only works in the interior, which is why we obtain
	\begin{align*}
		\int_{\Om_h^n} \nabla \mu_h^n \nabla(\eta^{\ell_n} - \eta_h) =&\ \sum_{T \in \calT_h^n} \int_T \nabla \mu_h^n \nabla(\eta^{\ell_n} - \eta_h) \nonumber \\
		 =&\ \sum_{S \in E_{\text{int}}^n} \int_S [\nabla \mu_h^n \cdot \nu] (\eta^{\ell_n} - \eta_h) + \sum_{S \in E_{\Ga_h}^n} \int_S \partial_\nu \mu_h^n (\eta^{\ell_n} - \eta_h),
	\end{align*}
	where the last term is matched with the remaining boundary part in $(\eta^{\ell_n} - \eta_h)$ to obtain an estimator on the boundary condition \eqref{eq:robintypecond2}:
	\begin{align*}
		\sum_{S \in E_{\Ga_h}^n} \int_S \partial_\nu \mu_h^n (\eta^{\ell_n} - \eta_h) - L^{-1} \int_{\Ga_h^n} (\beta\theta_h^n-\mu_h^n) (\eta^{\ell_n} - \eta_h) \nonumber \\
		=\sum_{S \in E_{\Ga_h}^n} \int_S \Big(\partial_\nu \mu_h^n- L^{-1} (\beta\theta_h^n-\mu_h^n)\Big) (\eta^{\ell_n} - \eta_h).
	\end{align*}
	Then we can bound the resulting terms using Cauchy-Schwarz and Hölder's inequality and finally using the interpolation estimates from \cite[Corollary~3.4]{CamachoDemlow2015} and norm equivalence (see Remark~\ref{rem:normequivalence}), to obtain the result:
	\begin{equation*}
		m(R_h^1(t),\zeta) \leq C \eta_{h,1}^n \| \zeta \|_{\calV}. \label{eq:Res1 - space est}
	\end{equation*}
	
	In order to estimate the geometric part of the residual we use Lemma~\ref{Lemma:GeomAppError} element-wise, to obtain
	\begin{align*}
		m(R_g^1(t),\zeta) \leq&\ C\Big( \sum_{(T,S)\in B_h^n \times E_{\Ga_h}^n} (h_T^2,h_S^4) \| \partial_{\tau_n} u_h^n \|_{m_h^{n,T,S}}^2 \nonumber \\
		&\ + \sum_{(T,S) \in B_h^n \times E_{\Ga_h}^n} (h_T^2,h_S^4) \| \nabla w_h^n \|_{m_h^{n,T,S}}^2 \\
		&\ + \sum_{S \in E_{\Ga_h}^n} h_S^4 \| \beta \theta_h^n - \mu_h^n \|_{L^2(S)}^2 \Big)^{\frac12} \| \zeta \|_{\calV}.
	\end{align*}
	One now easily sees, that the resulting terms match to the already obtained spatial estimators, with the exception of the jump terms of the spatial estimator, where we are not able to obtain the analogous term which stems from an additional partial integration. This yields
	\begin{align}
		m(R_g^1(t),\zeta) \leq C (\eta_{h,1}^n + \calG_{g,1}^n) \| \zeta \|_{\calV}. \label{eq:Res1 - geom est}
	\end{align}
	
	Finally the coarsening residual $R_c^1$ is estimated by splitting the error in $(u_h^{n-1})^{\ell_{n-1}}$ and $(\calI_{n-1}^n u_h^{n-1})^{\ell_n}$ into 
	\begin{align*}
		(u_h^{n-1})^{\ell_{n-1}} - (\calI_{n-1}^n u_h^{n-1})^{\ell_n} =&\ \bigg[ \Big((u_h^{n-1})^{\ell_{n-1}} - (\calI_{n-1}^\oplus u_h^{n-1})^{\ell_\oplus}\Big) \\
		&\  + \Big( (\calI_{n-1}^\oplus u_h^{n-1})^{\ell_\oplus} - (\calI_{n}^\oplus \calI_{n-1}^n u_h^{n-1})^{\ell_\oplus}\Big) \\
		&\ + \Big( (\calI_{n}^\oplus \calI_{n-1}^n u_h^{n-1})^{\ell_\oplus} - (\calI_{n-1}^n u_h^{n-1})^{\ell_n}\Big)\bigg],
	\end{align*}
	where the first and third difference are of geometric nature, as the underlying finite element functions are identical, yet lifted differently, and the error terms according to Lemma~\ref{Lemma:common_interpolation} are collected in $\calG_{c,1}^n$. The second term is the contribution of the coarsening estimator $\eta_{c,1}^n$ and captures the information loss in $\calI_{n-1}^n$ by coarsening the mesh from one step to the next. This leads to
	\begin{align}
		m(R_c^1(t),\zeta) \leq C (\eta_{c,1}^n + \calG_{c,1}^n) \| \zeta \|_{\calV}. \label{eq:Res1 - coarse est}
	\end{align}
	Combining \eqref{eq:Res1 - space est}, \eqref{eq:Res1 - geom est} and \eqref{eq:Res1 - coarse est} gives the estimate for the dual norm.
	
	The estimate for the mass of the residual is obtained in an analogous way: We split the expression as in \eqref{eq:Res 1 split}, with $\zeta=(\beta,1)$, to arrive at
	\begin{align*}
		m(R^1(t),(\beta,1)) =&\ \frac{1}{\tau_n} m\big((\calI_{n-1}^n u_h^{n-1})^{\ell_n} - (u_h^{n-1})^{\ell_{n-1}},(\beta,1)\big) \nonumber \\
		&\ + \big[m((\partial_{\tau_n} u_h^n)^{\ell_n}, (\beta,1)) - m_h^n(\partial_{\tau_n} u_h^n, (\beta,1)^{-\ell_n})\big],
	\end{align*}
	where we have used the fact that $(\beta,1)^{-\ell_n} \in \calV_h^n$ to eliminate the spatial residual. Then this is simply estimated by the previously developed techniques for the coarsening and geometric residual to the proposed bound. Note that the contributions to the geometric indicator $\calG_{g,1}^n$ can be again crudely estimated by $\eta_{h,1}^n$ and $\eta_{c,1}^n$.
\end{proof}

\begin{rem}
	\begin{enumerate}
		\item Note that the estimate of the geometric residual by the spatial estimator is crude, as for the boundary terms we neglected an additional $h$ and the sum for the bulk terms only contains elements which have at least two nodes on the boundary. 
		\item We observe that by the definition of the common mesh and the respective interpolation operators, the difference 
		\begin{equation*}
			\calI_{n-1}^\oplus u_h^{n-1} - \calI_{n}^\oplus \calI_{n-1}^n u_h^{n-1},
		\end{equation*}
		vanishes exactly, when there is no coarsening from $\Om_h^{n-1}$ to $\Om_h^n$, cf. Lemma~\ref{Lemma:common_interpolation}.
	\end{enumerate}
\end{rem}

The second residual in \eqref{eq:DefResidualNorms} is estimated by the exact same techniques, with some slight changes, resulting in the analogous estimators 
\begin{subequations} \label{eq:Error indicators 2}
	\begin{align}
		\eta_2^n :=&\ \eta_{t,2}^n + \eta_{F,2}^n + \eta_{h,2}^n + \eta_{c,2}^n,\\
		(\eta_{t,2}^n)^2 :=&\ \sum_{(T,S) \in \calT_h^n \times E_{\Ga_h}^n} \| \eps \nabla (u_h^n - \calI_{n-1}^n u_h^{n-1})\|_{m_h^{n,T,S}}^2	+ \|\inveps (u_h^n - \calI_{n-1}^n u_h^{n-1}) \|_{m_h^{n,T,S}}^2 \nonumber \\
		&\ \hphantom{\sum_{(T,S) \in \calT_h^n \times E_{\Ga_h}^n}} + \| K^{-1}(\alpha \psi_h^{n-1} - \phi_h^{n-1}) - K^{-1}(\alpha \psi_h^{n} - \phi_h^{n}) \|_{m_h^{n,T,S}}^2,, \\
		(\eta_{F,2}^n)^2 :=&\ \sum_{(T,S) \in \calT_h^n \times E_{\Ga_h}^n} h^2 \| \inveps \nabla \calI_{n-1}^n F'(u_h^{n-1}) \|_{m_h^{n,T,S}}^2  \\
		(\eta_{h,2}^n)^2 :=&\ \sum_{(T,S) \in \calT_h^n \times E_{\Ga_h}^n} h^2 \| w_h^n - \inveps \calI_{n-1}^n F'(u_h^{n-1}) \|_{m_h^{n,T,S}}^2 \nonumber \\
		&\ + \sum_{(S,P) \in E_{\text{int}}^n \times P_{\Ga_h}^n} h \| \eps [\nabla u_h^n \cdot \nu]\|_{m_h^{n,S,P}}^2 \nonumber \\
		&\ + \sum_{S \in E_{\Ga_h}^n } h_S \| \eps\partial_\nu \phi_h^n - K^{-1} (\alpha \psi_h^n - \phi_h^n)\|_{L^2(S)}^2 + h_S^2 \| K^{-1}(\alpha \psi_h^n -  \phi_h^n) \|_{L^2(S)}^2  \\
		(\eta_{c,2}^n)^2 :=&\  \sum_{(T,S) \in \calT_h^\oplus \times E_{\Ga_h}^\oplus} \|\eps \nabla \big( \calI_{n-1}^\oplus u_h^{n-1} - \calI_{n}^\oplus \calI_{n-1}^n u_h^{n-1}\big)\|_{m_h^{n,T,S}}^2  \nonumber  \\
		&\ \hphantom{\sum_{(T,S) \in \calT_h^\oplus \times E_{\Ga_h}^\oplus}}  + \| \inveps (\calI_{n-1}^\oplus u_h^{n-1} - \calI_{n}^\oplus \calI_{n-1}^n u_h^{n-1})\|_{m_h^{n,T,S}}^2 \nonumber   \\
		&\ \hphantom{\sum_{(T,S) \in \calT_h^\oplus \times E_{\Ga_h}^\oplus}}  +\| K^{-1} \big(\alpha(\calI_{n-1}^\oplus \psi_h^{n-1} - \calI_{n}^\oplus \calI_{n-1}^n \psi_h^{n-1})  \nonumber \\ 
		&\ \hphantom{\sum_{(T,S) \in \calT_h^\oplus \times E_{\Ga_h}^\oplus}} \qquad  - (\calI_{n-1}^\oplus \phi_h^{n-1} - \calI_{n}^\oplus \calI_{n-1}^n \phi_h^{n-1})\big)\|_{m_h^{n,T,S}}^2,
	\end{align}
	for $n=1,\dots,N$, with geometric estimators given by
	\begin{align}
		\calG_2^n :=&\ \calG_{g,2}^n + \calG_{c,2}^n \\
		(\calG_{g,2}^n)^2 :=&\ \sum_{(T,S) \in B_h^n \times E_{\Ga_h}^n} (h_T^2,h_S^4) \| \eps \nabla u_h^n \|_{m_h^{n,T,S}}^2 \\
		(\calG_{c,2}^n)^2:=&\ \sum_{(T,S) \in \hat B_h^{n-1} \times \hat E_h^{n-1}}   (\eps^2 h^2 + \eps^{-2} h^4) \| \nabla u_h^{n-1} \|_{m_h^{n-1 ,T,S}}^2 \nonumber \\
		&\ + \sum_{S \in \hat E_{\Ga_h}^{n-1}}  h_S^4 \| K^{-1} \nabla_{\Ga_h} (\alpha\psi_h^{n-1} - \phi_h^{n-1}) \|_{L^2(S)}^2  \nonumber \\
		&\ + \sum_{(T,S) \in \hat B_h^{n} \times \hat E_h^{n}}   (\eps^2 h^2 + \eps^{-2} h^4) \| \nabla (\calI_{n-1}^n u_h^{n-1}) \|_{m_h^{n,T,S}}^2  \nonumber \\
		&\ + \sum_{S \in \hat E_{\Ga_h}^{n}}  h_S^4 \| K^{-1}  \nabla_{\Ga_h} (\alpha \calI_{n-1}^n \psi_h^{n-1} - \calI_{n-1}^n \phi_h^{n-1}) \|_{L^2(S)}^2. 
	\end{align}
\end{subequations}

\begin{lem} \label{Lemma:Residual bound 2}
	Let $n=1,\dots,N$ and $t \in (t^{n-1},t^n]$, then the dual norm \eqref{eq:dualNormDef} of the second residual $R^2$ fulfills
	\begin{align}
		\| R^2(t) \|_{\calV^{\ast}} \leq C \Big(\eta_2^n + \calG_2^n\Big)
	\end{align}
	for some constant $C>0$ independent of $h$, $\tau$ and $\eps$, and the error indicators defined in \eqref{eq:Error indicators 2}.
\end{lem}

\begin{proof}
	As in the proof of Lemma~\ref{Lemma:Residual bound 1}, we split the expression \eqref{eq:DefResidual2} into
	\begin{align*}
		m(R^2,\eta) =&\ m(w_{h,\tau},\eta) - a^{K,\alpha,\eps}(u_{h,\tau},\eta) - m(\inveps F'(u_{h,\tau}) ,\eta) \nonumber \\
		=&\ m(R_c^2,\zeta) + m(R_g^2,\zeta) + m(R_t^2,\zeta) + m(R_h^2,\zeta).
	\end{align*}
	with the coarsening, geometric, temporal and spatial residuals
	\begin{align*}
		m(R_c^2,\zeta):=&\ - a^{K,\alpha,\eps}(u_{h,\tau} - \widetilde u_{h,\tau}^{\ell_n},\eta) - m(\inveps F'(u_{h,\tau}) - \inveps F'(\widetilde u_{h,\tau}^{\ell_n}),\eta), \\
		m(R_g^2,\zeta):=&\ \Big[m((w_h^n)^{\ell_n},\eta) - m_h^n( w_h^n,\eta^{-\ell_n})\Big] - \Big[ a^{K,\alpha,\eps}(\widetilde u_{h,\tau}^{\ell_n},\eta) - a_h^{n,K,\alpha,\eps}(\widetilde u_{h,\tau},\eta^{-\ell_n}) \Big] \nonumber \\
		&\ - \Big[ m(\inveps F'(\widetilde u_{h,\tau}^{\ell_n}),\eta) - m_h^n(\inveps F'(\widetilde u_{h,\tau}),\eta^{-\ell_n})\Big], \\
		m(R_t^2,\zeta):=&\ - a_h^{n,K,\alpha,\eps}(\widetilde u_{h,\tau} - u_h^n,\eta^{-\ell_n}) - m_h^n(\inveps F'(\widetilde u_{h,\tau}) - \inveps \calI_{n-1}^n F'(u_h^{n-1}),\eta^{-\ell_n}), \\
		m(R_h^2,\zeta):=&\  m_h^n(w_h^n,\eta^{-\ell_n}) - a_h^{n,K,\alpha,\eps}(u_h^n,\eta^{-\ell_n}) - m_h^n(\inveps \calI_{n-1}^n F'(u_h^{n-1}),\eta^{-\ell_n}).
	\end{align*}
	The terms $R_c^2$, $R_g^2$ and $R_t^2$ are then estimated exactly as in Lemma~\ref{Lemma:Residual bound 1}, with the only notable difference being the non-linearity, which can be treated by the proposed Lipschitz continuity. 
	
	The first term in the newly appearing temporal estimator is simply estimated by the maximum of the linear interpolation in time, which gives the first part of the estimator $\eta_{t,2}^n$. The only substantially new term appears in the second term of $R_t^2$, coming from the additional interpolation operator $\calI_{n-1}^n$ in front of $F'$, and is estimated element-wise into the new estimator $\eta_{F,2}^n$, 
	\begin{align*}
		&m_h(\inveps F'(\widetilde u_{h,\tau}) - \inveps \calI_{n-1}^n F'(u_h^{n-1}),\eta^{-\ell}) = \nonumber \\
		  &\inveps m_h(F'(\widetilde u_{h,\tau}) - F'(\calI_{n-1}^n u_h^{n-1}),\eta^{-\ell})  + \inveps m_h(F'(\calI_{n-1}^n u_h^{n-1}) - \calI_{n-1}^n F'(\calI_{n-1}^n u_h^{n-1}),\eta^{-\ell}) \\
		&\ \leq C \inveps \Big(\| u_h^n - \calI_{n-1}^n u_h^{n-1}\|_{m_h} + h \| \nabla \calI_{n-1}^n F'(\calI_{n-1}^n u_h^{n-1}) \|_{m_h}\Big) \| \eta \|_{\calV},
	\end{align*}
	by a standard interpolation estimate for $\calI_{n-1}^n$, see, e.g. \cite{BrennerScott2008}, and the proposed Lipschitz continuity of $F'$. The resulting estimators are collected in \eqref{eq:Error indicators 2} and the analogous stated bound for $R^2$ is obtained.
\end{proof}

\subsection{Main result: A posteriori error estimates \label{sec:APosterioriErrorEstimates}}
Using the bounds of Lemma~\ref{Lemma:Residual bound 1} and Lemma~\ref{Lemma:Residual bound 2}, we get the full residual bounds in terms of the error indicators: 
\begin{align}
	\eps^{\frac12} \|R^1\|_{L^2(0,t^n; \calV^{\ast})} &\leq C \eps^{\frac12} \sqrt{t^n} \Big(\sup_{0\leq k \leq n} \eta_{1}^k + \sup_{0\leq k \leq n} \mathcal{G}_1^k \Big), \\
	\eps^{-\frac12} \| m(R^1,(\beta,1))\|_{L^2(0,T)} &\leq C \eps^{-\frac12} \sqrt{t^n} \Big(\sup_{0\leq k \leq n} \eta_{M,1}^k + \sup_{0\leq k \leq n} \mathcal{G}_{1}^k \Big), \\
	\eps^{-\frac12} \|R^2\|_{L^2(0,t^n; \calV^{\ast})} &\leq C \eps^{-\frac12} \sqrt{t^n}  \Big(\sup_{0\leq k \leq n} \eta_{2}^k + \sup_{0\leq k \leq n} \mathcal{G}_2^k \Big). 
\end{align}

By rearranging the estimators we relate them to their temporal, spatial, coarsening and geometric components,
\begin{equation}\label{eq:Final_Estimators}
	\begin{aligned}
		\eta^n:=&\  \eta_t^n + \eta_h^n + \eta_c^n, \\
		(\eta_t^n)^2 :=&\ \sup_{0\leq k \leq n} \inveps t^n \eta_{t,2}^k, \\
		(\eta_h^n)^2 :=&\ \sup_{0\leq k \leq n} \eps \, t^n \eta_{h,1}^k + \sup_{0\leq k \leq n} \eps^{-1}  t^n \big(\eta_{h,2}^k + \eta_{F,2}^k\big) , \\
		(\eta_c^n)^2 :=&\ \sup_{0\leq k \leq n} \eps \,  t^n \eta_{c,1}^k + \sup_{0\leq k \leq n} \eps^{-1} t^n \big(\eta_{M,1}^k + \eta_{c,2}^k\big), \\
		(\calG^n)^2 :=&\ \sup_{0\leq k \leq n} \eps^{-1} t^n \big(\calG_{1}^k + \calG_{2}^k\big).
	\end{aligned} 
\end{equation}

Finally we use \eqref{eq:ErrorByResidualBound} of Lemma~\ref{Lemma:error_residual_bound} to infer the final a posteriori error estimate, with the initial error given by
\begin{align}
\label{eq:initial_error}
	\calE_I := | u^0-(u_h^0)^{\ell_0} |_{\calH_{L,\beta}^{-1}}^2 + \eps \| u^0-(u_h^0)^{\ell_0} \|_{\calV}^2.
\end{align}

\begin{thm}\label{thm:MainResultReliability}
	Let $(u,\psi,\mu,\theta)$ be a solution of the bulk--surface Cahn--Hilliard equation with dynamic boundary conditions \eqref{eq:CHreact}, with nonlinear potentials satisfying \eqref{eq:RegularityF}, then there exists $h_0>0$, such that for all $h\leq h_0$, the linear interpolation $(u_{h,\tau},\psi_{h,\tau},\mu_{h,\tau},\theta_{h,\tau})$ of the solution of the numerical scheme \eqref{eq:CHweakD} satisfies the following
	error bound, containing the estimators \eqref{eq:Final_Estimators}, for $0<t^n=\tau_1+\cdots+\tau_n\leq T$:
	\begin{align*}
		\|u-u_{h,\tau}\|_{L^\infty(0,t^n;\calH_{L,\beta}^{-1})} + \eps \|u- u_{h,\tau}\|_{L^2(0,t^n;\calV)} \leq c^\ast \bigg(\eta^n + \mathcal{G}^n + \calE_I \bigg)  ,
	\end{align*}
	where $c^\ast > 0$ is independent of $h$ and $\tau_1,\dots,\tau_n$, but depends on norms of the exact solution, parameters of the problem \eqref{eq:CHreact} and the domain $\Omega$.
\end{thm}

\begin{rem}
	Note that, the exponential dependence of the constant in Theorem~\ref{thm:MainResultReliability} essentially limits its use in practice to time scales of $\mathcal{O}(\eps^3)$. As noted in the introduction, this common problem for phase field models was overcome in \cite{FengWu2007,BartelsMuller2010,BartelsMullerOrtner2011} for the Cahn--Hilliard equations under the condition that the interface only undergoes certain topological changes. These arguments depend heavily on spectral estimates for a linearized Cahn--Hilliard operator, which are currently not available for the present case of dynamic boundary conditions. However, from our numerical experiments we expect that the dependence on the parameter $\eps$ could be improved.
\end{rem}

\section{Numerical Experiments} \label{sec:Numerical Experiments}
We want to showcase the potential of using adaptivity when approximating the Cahn-Hilliard problem \eqref{eq:CHreact}.

In the upcoming numerical experiments, we computed the adaptive solution based on the error indicators derived in Section~\ref{sec:indicators_bound_residual_bound}.

We use a minimal setup for adaptivity, which consists of the following steps
\begin{enumerate}
\item Determine an initial coarse base mesh, such that the assumptions in Section~\ref{sec:discretization} are fulfilled.
\item Refine the mesh for $t^0$ such that the initial error, given in \eqref{eq:initial_error}, is bounded by $\mathcal{E}_I \leq \TOL$.
\item For each time step
\begin{enumerate}
\item Solve the fully discrete problem \eqref{eq:CHweakD1}, \eqref{eq:CHweakD2} and compute the error indicators \eqref{eq:Final_Estimators}.
\item While the temporal error is too large, i.e.\ $\eta_t  > \TOL_\tau$, reduce time step size.
\item While the spatial error is too large, i.e.\ $\eta_h > \TOL_h$, refine spatially.
\end{enumerate}
\item Store solution and coarsen time, if the temporal error is sufficiently small (we checked $\eta_\tau < 0.9\TOL_\tau$), and space (as long as $\eta_c < \TOL_c$) for next step.
\end{enumerate}

We note that we simplified the numerical computations by omitting some contributions of \eqref{eq:Final_Estimators} we will briefly comment on those:
\begin{enumerate}
\item For the initial error we assume (and use) an exact solution. Thus step 2 is fulfilled. We further follow that within the first step of solving the fully discrete problem, we update the required $\mathcal{I}_{0}^1 u_h$ based on the exact initial data, instead of Lagrangian interpolation.  
\item The geometric indicators $\mathcal{G}^n$ of \eqref{eq:Final_Estimators}, divide in two parts. Whereas both are generally high-order in space, the mesh-transfer parts include inverse powers of the temporal stepsize $\tau$.
The typical geometric perturbations arising by the shift from discrete to continuous bilinear form, namely $\mathcal{G}_{g,1}^n$ and $\mathcal{G}_{g,2}^n$, are asymptotically neglectable.
The first contribution of each term is justified to be small by the fact that for $C^2$-domains one can observe that the geometric contributions are asymptotically high-order, see for example \cite{DemlowDziuk2007} for surface terms. And the additional contributions are $h$-scaled versions of terms in other indicators.
Finally the mesh-transfer errors, represented by $\mathcal{G}_{c,1}^n$ and $\mathcal{G}_{c,2}^n$ measure mesh-changes and thus hinder strong coarsening and refining on the boundary. Whereas the second is easy to manage, the first can require reverting to previous meshes due to the $\tau_n^{-2}$ factor, which is algorithmically manageable, but expensive. See the discussion in \cite[Section~6]{Lantelme2026} for more details. 
Herein we did not implement a routine guaranteeing the mesh-transfer to be under a specified tolerance, but checked heuristically that they are non-dominant. 
\end{enumerate}

The upcoming experiments include: In Section~\ref{sec:ex1} illustrative examples, showcasing the functionality of coarsening and refinement of the interface, evolution of mass and evolution of stepsize.
In Section~\ref{sec:ex2} we provide a reliability and efficiency plot, showcasing not only the proven upper bound \eqref{thm:MainResultReliability} but also heuristically suggesting a lower bound.
Finally, in Section~\ref{sec:ex3} we performed a runtime comparison, benchmarking the adaptive routine against a standard uniform routine.

\subsection{Evolution experiments}
\label{sec:ex1}
In this section we present evolution experiments on a curved domain $\Om$, which is given as a perturbed circle via the signed distance function in polar coordinates $(r,\phi) \in (0,\infty) \times [0,2\pi)$,
\begin{align*}
	d(r,\phi)= r - (1+0.04\cdot \cos(5\phi)).
\end{align*} 
We consider the following situations:
\begin{enumerate}
	\item The initial data 
	$$u^0(x)=\sin(4\pi x_1) \cos(4\pi x_2),$$
	taken from \cite[Section~7.2.2]{KnopfLam2020}, and the parameters $K=10^{-5}$, $L=10^{-5}$ and $\beta=1$.
	\item The initial data
	$$u^0(x) = \tanh(((x_1-0.3)^2 + x_2^2-0.25^2)/\eps) \cdot \tanh(((x_1 + 0.3)^2 + x_2^2-0.3^2)/\eps),$$
	taken from \cite[Section~6]{FengWu2007}, and the parameters $K=10^{-5}$, $L=10^{5}$ and $\beta=1$.
\end{enumerate}
In both experiments we choose the diffuse interface parameter $\eps=0.02$, a starting mesh with $103$ nodes and a maximal time step size $\tau_{\max}=4\cdot 10^{-5}$. 

\begin{figure}[!t]
	\centering\includegraphics[width=\textwidth]{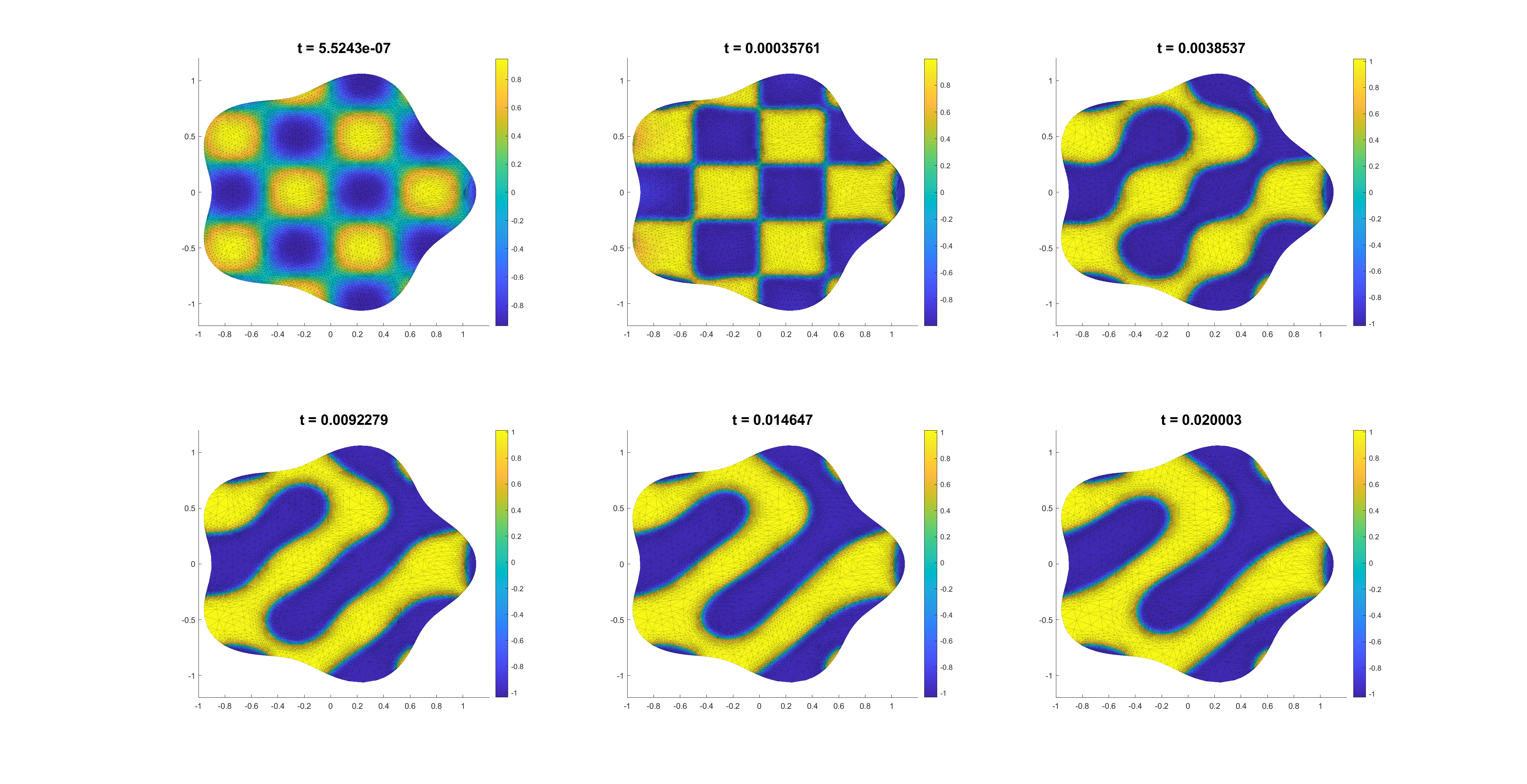}
	\caption{Evolution plot for the bulk phase field $\phi$ of the numerical solution of our adaptive algorithm for the initial data 1).}
	\label{fig:EvolutionPlot - Sin Cos}
\end{figure}

In the first experiment we have chosen the tolerances to be 
$$[\TOL_\tau,\TOL_h,\TOL_c] = [60, 1, 1],$$
and Figure~\ref{fig:EvolutionPlot - Sin Cos} shows snapshots of the evolution of the bulk phase field $\phi$. We observe that the adaptive mesh follows nicely the interface between the two phases, and the solution shows the expected evolution behavior, see, e.g., \cite{KnopfLam2020,BullerjahnKovacs2024} for similar experiments on uniform approximations. 

\begin{figure}[!t]
	\centering\includegraphics[width=\textwidth]{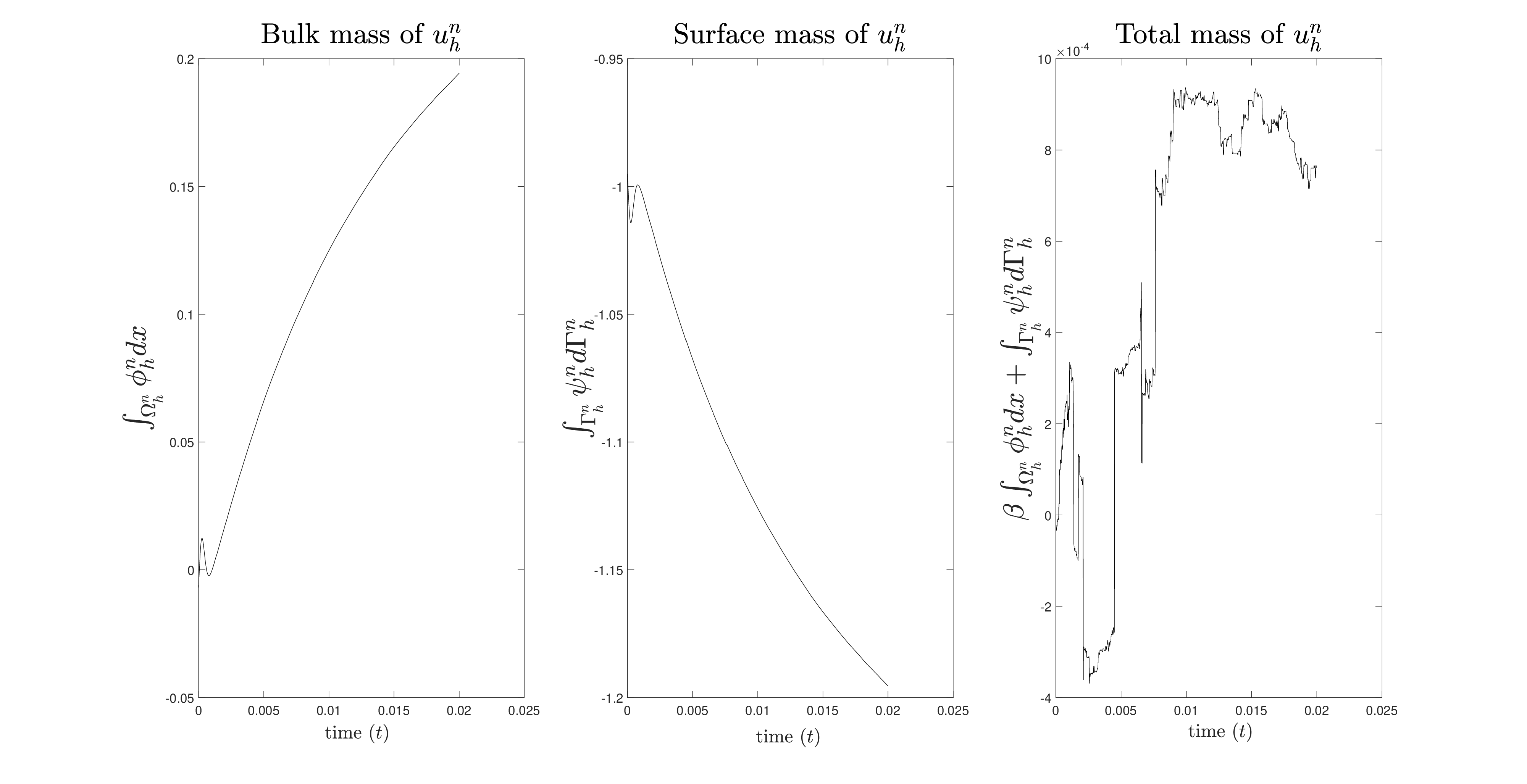}
	\caption{Evolution of the bulk, surface and combined mass of $u$ for the numerical solution of our adaptive algorithm for the initial data 1) (Note the axis scaling).}
	\label{fig:MassPlot - Sin Cos}
\end{figure}

In Figure~\ref{fig:MassPlot - Sin Cos} we report on the change in mass of the solution $u=(\phi,\psi)$. Note that the exact solution satisfies a mass conservation of the combined bulk--surface mass (see \cite{KnopfStange2024}), which is no longer exactly preserved in our adaptive scheme, cf. Remark~\ref{rem:mass conservation}, in contrast to the uniform approximation, see, e.g. \cite[Section~11.2]{Bullerjahn2025}.

\begin{figure}[!t]
	\centering\includegraphics[width=\textwidth]{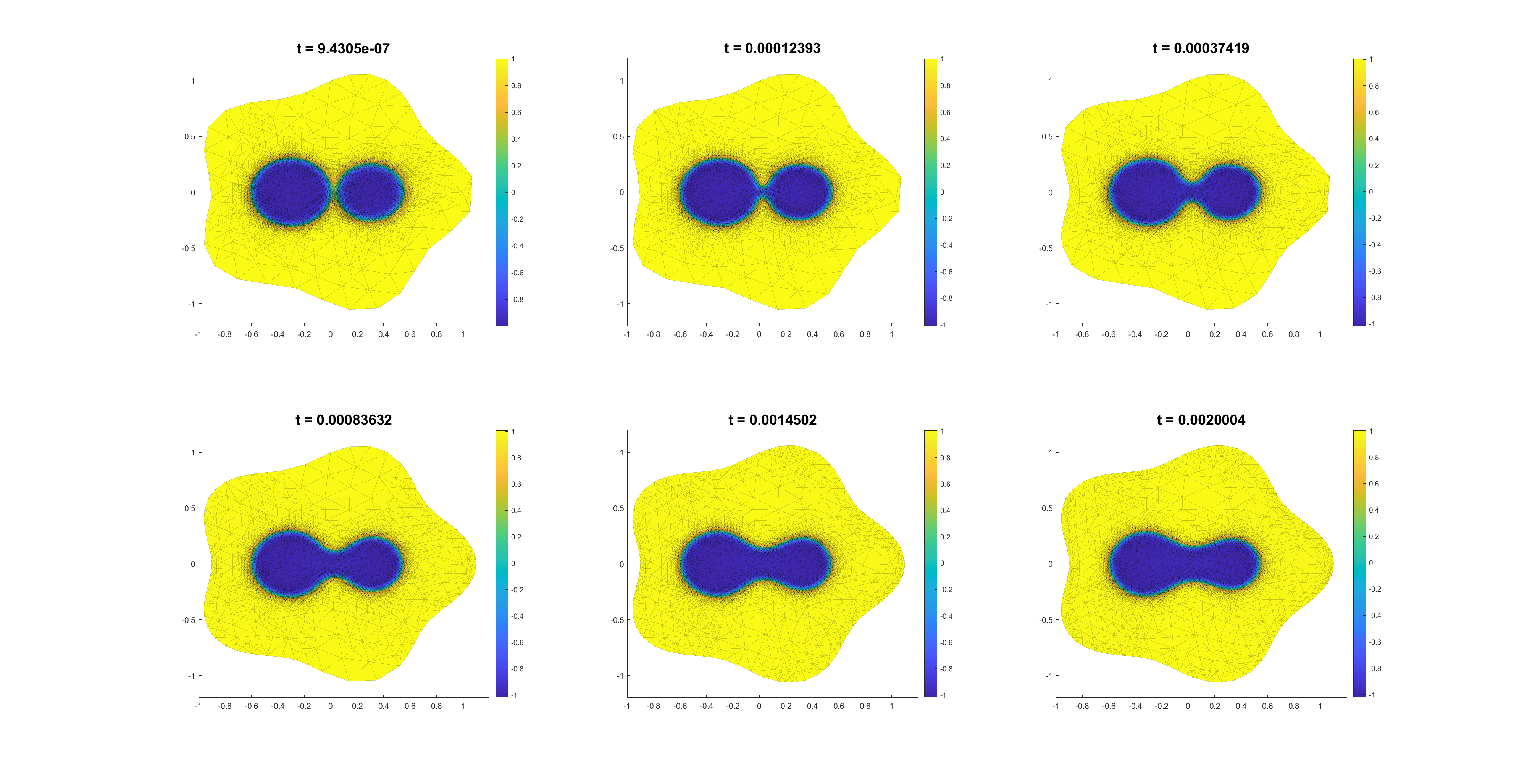}
	\caption{Evolution plot for bulk phase field $\phi$ of the numerical solution of our adaptive algorithm for the initial data 2).}
	\label{fig:EvolutionPlot - MergeCircles U}
\end{figure}

\begin{figure}[!t]
	\centering\includegraphics[width=\textwidth]{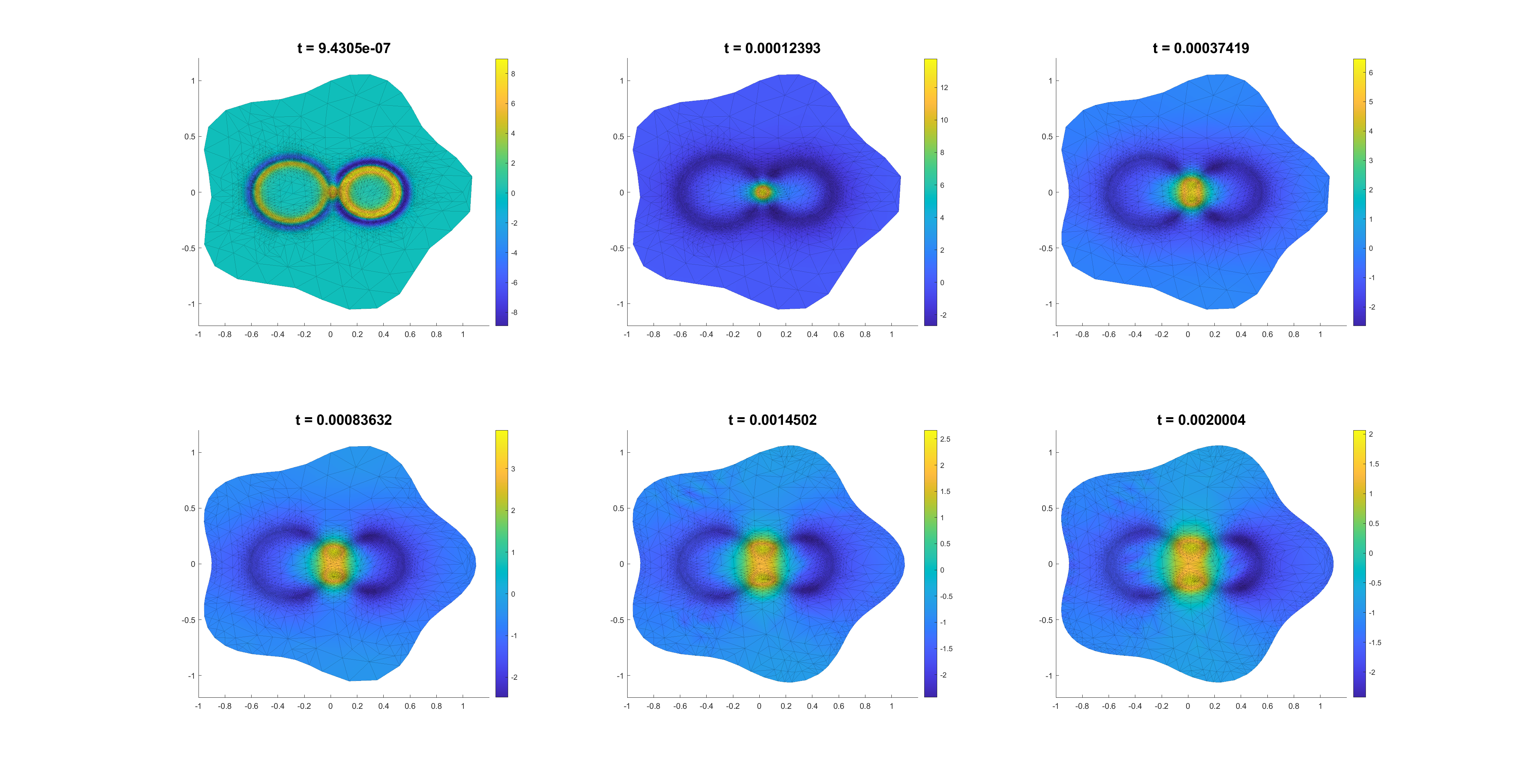}
	\caption{Evolution plot for bulk chemical potential $\mu$ of the numerical solution of our adaptive algorithm for the initial data 2).}
	\label{fig:EvolutionPlot - MergeCircles W}
\end{figure}

For the second experiment we have chosen the tolerances to be 
$$[\TOL_\tau,\TOL_h,\TOL_c] = [30, 1, 1],$$
and Figures~\ref{fig:EvolutionPlot - MergeCircles U} \& \ref{fig:EvolutionPlot - MergeCircles W} show snapshots of the evolution of the bulk phase field $\phi$, respectively the bulk chemical potential $\mu$. We observe that the adaptive mesh follows nicely the interface between the two phases of $\phi$, while also resolving $\mu$, if necessary. The solution shows the expected evolution behavior, see, e.g., \cite{FengWu2007} for similar experiments.

\begin{figure}[!t]
	\centering\includegraphics[width=\textwidth]{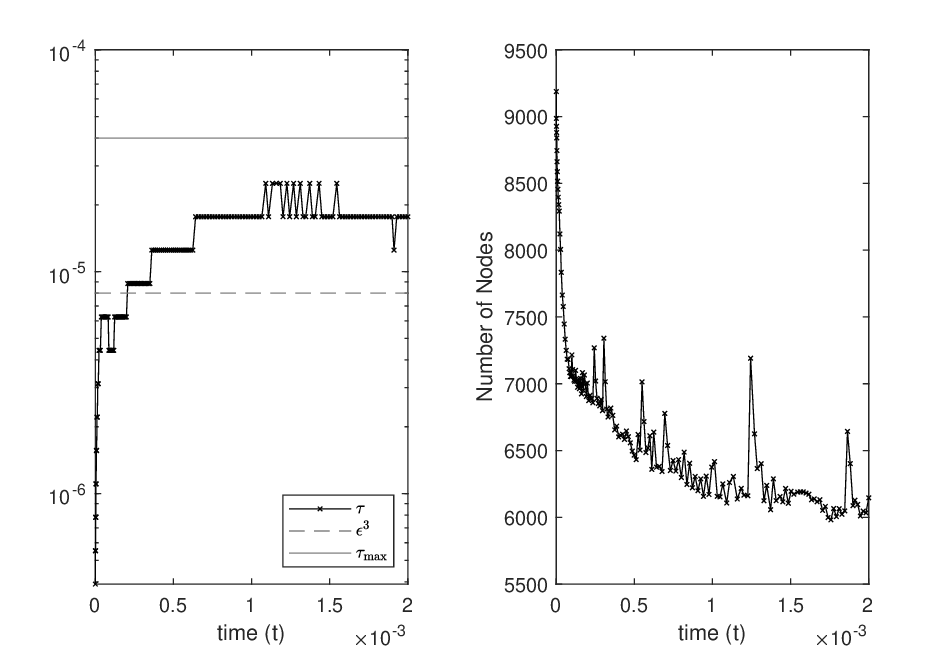}
	\caption{Evolution of the time step size $\tau$ and the number of nodes $N$ for the numerical solution of our adaptive algorithm for the initial data 2).}
	\label{fig:NTauPlot - MergeCircles}
\end{figure}

In Figure~\ref{fig:NTauPlot - MergeCircles} we report on the evolution of the time step size $\tau$ and the number of nodes $N$ over time. We observe an exponential decay of both quantities followed by a saturation, which is expected due to the initial regularization on a fast time scale, followed by a coarsening behavior on a slower time scale, cf. Remark~\ref{rem:Cont results Epsilon}. 

\subsection{Reliability and efficiency}
\label{sec:ex2}
Whereas our main result \ref{thm:MainResultReliability} only proves reliability of the error indicator (up to oscillation) we can heuristically analyse both reliability and efficiency. 

We computed the error indicators, and the error with respect to a reference solution for a set of (almost) uniform spatial meshes and fixed time steps. 
The problem was posed on the unit disk and we analysed the droplet example presented in \cite[Section~6]{KnopfLamLiuMetzger2021} with the following parameter choices:
The elliptic droplet is centered at the boundary in $(x^0,y^0) =(-1,0)$, its initial horizontal and vertical width is $a_1 =0.75$ and $a_2 =0.18$ respectively, and interface parameter $\delta = 0.5$.  
Thus the initial condition is written as
\begin{align}
u^0(x) = \tanh\Bigg(\frac{-1}{\delta} \bigg(\frac{\big(x_1-x^0\big)^2}{a_1^2}+\frac{\big(x_2-y^0\big)^2}{a_2^2}-1\bigg) \Bigg). \label{eq:Num initial data droplet}
\end{align}

Additionally, we set the parameters $\eps = 0.1$, $K = L = 1$ and $T = 0.01$.
We determined the dominant parts of the error estimators \eqref{eq:Final_Estimators} for a set of uniform meshes with maximal mesh-width $h$ and time step sizes $\tau$:
\begin{align*}
h &\in \{ 0.4300,    0.3487,    0.1807,    0.0918,    0.0463,    0.0232\}, \\
\tau &\in \{ 10^{-2}, 10^{-3}, 10^{-4}\}.
\end{align*}

To check the reliability we approximate the $L^\infty(\mathcal{H}^{-1}_{L,\beta})$ and the $L^2(H^1)$-error of the bulk--surface phase field for each mesh and time step pair.
This is done by comparing the approximation on a fine reference mesh to a reference solution.
The reference solution was computed with the step-size $\tau = 10^{-6}$ and an almost uniform mesh with $123457$ vertices ($h =0.0116$).

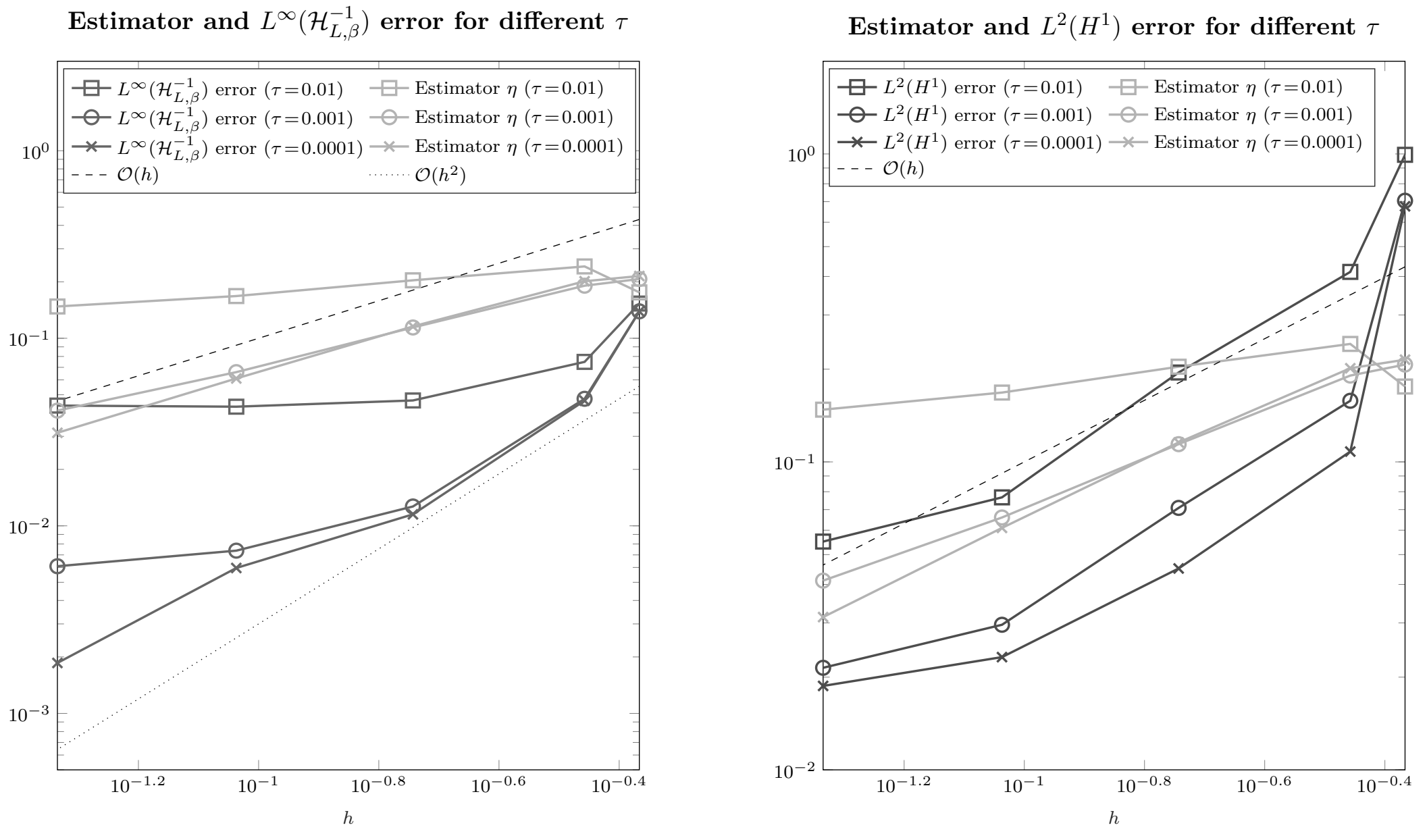
\begin{figure}[htbp]
	\centering
	\resizebox{\textwidth}{!}{%
%
%
\begin{tikzpicture}

\begin{axis}[%
width=3.486in,
height=4.245in,
at={(1.354in,0.573in)},
scale only axis,
xmode=log,
xmin=0.0462624701308489,
xmax=0.430042983911045,
xminorticks=true,
xlabel style={font=\color{white!15!black}},
xlabel={$h$},
ymode=log,
ymin=0.0005,
ymax=3,
yminorticks=true,
axis background/.style={fill=white},
title style={font=\bfseries},
title={\large{Estimator and $L^\infty(\mathcal{H}_{L,\beta}^{-1})$  error for different $\tau$}},
legend style={at={(0.01,0.99)}, anchor=north west, legend cell align=left, align=left, draw=white!15!black},
legend columns=2,
]
\addplot [color=white!40!black, mark=square, mark options={solid, white!40!black},line width=1pt,mark size=3pt]
  table[row sep=crcr]{%
0.430042983911045	0.153159524609191\\
0.348693200988337	0.0748537821022449\\
0.180653453374279	0.0466109408560977\\
0.0918101812671815	0.0431794511571898\\
0.0462624701308489	0.0437687873732793\\
};
\addlegendentry{$L^\infty(\mathcal{H}_{L,\beta}^{-1})$ error $(\tau \!=\! 0.01)$}

\addplot [color=white!70!black, mark=square, mark options={solid, white!70!black},line width=1pt,mark size=3pt]
  table[row sep=crcr]{%
0.430042983911045	0.175695234977988\\
0.348693200988337	0.241583793078715\\
0.180653453374279	0.203741273948975\\
0.0918101812671815	0.167933967743723\\
0.0462624701308489	0.147768955099406\\
};
\addlegendentry{Estimator $\eta$ $(\tau \!=\! 0.01)$}

\addplot [color=white!40!black, mark=o, mark options={solid, white!40!black},line width=1pt,mark size=3pt]
  table[row sep=crcr]{%
0.430042983911045	0.1391951290364\\
0.348693200988337	0.0475736289669187\\
0.180653453374279	0.0126824146808769\\
0.0918101812671815	0.00737473865196632\\
0.0462624701308489	0.00608999009161801\\
};
\addlegendentry{$L^\infty(\mathcal{H}_{L,\beta}^{-1})$ error $(\tau \!=\! 0.001)$}

\addplot [color=white!70!black, mark=o, mark options={solid, white!70!black},line width=1pt,mark size=3pt]
  table[row sep=crcr]{%
0.430042983911045	0.207021323303769\\
0.348693200988337	0.19071651798131\\
0.180653453374279	0.114170704157295\\
0.0918101812671815	0.0659934381888048\\
0.0462624701308489	0.0411577757782331\\
};
\addlegendentry{Estimator $\eta$ $(\tau \!=\! 0.001)$}

\addplot [color=white!40!black, mark=x, mark options={solid, white!40!black},line width=1pt,mark size=3pt]
  table[row sep=crcr]{%
0.430042983911045	0.139046284668329\\
0.348693200988337	0.0464844126538016\\
0.180653453374279	0.0115196310061897\\
0.0918101812671815	0.00596000405310616\\
0.0462624701308489	0.00185761629732407\\
};
\addlegendentry{$L^\infty(\mathcal{H}_{L,\beta}^{-1})$ error $(\tau \!=\! 0.0001)$}

\addplot [color=white!70!black, mark=x, mark options={solid, white!70!black},line width=1pt,mark size=3pt]
  table[row sep=crcr]{%
0.430042983911045	0.214714895341424\\
0.348693200988337	0.20141392063202\\
0.180653453374279	0.115661015401895\\
0.0918101812671815	0.0610460721588877\\
0.0462624701308489	0.031337688478144\\
};
\addlegendentry{Estimator $\eta$ $(\tau \!=\! 0.0001)$}

\addplot [color=black, dashed]
  table[row sep=crcr]{%
0.430042983911045	0.430042983911045\\
0.348693200988337	0.348693200988337\\
0.180653453374279	0.180653453374279\\
0.0918101812671815	0.0918101812671815\\
0.0462624701308489	0.0462624701308489\\
};
\addlegendentry{$\mathcal{O}(h)$}

\addplot [color=black, dotted]
  table[row sep=crcr]{%
0.430042983911045	0.0554810904033345\\
0.348693200988337	0.0364760845246478\\
0.180653453374279	0.00979070106481579\\
0.0918101812671815	0.00252873281529382\\
0.0462624701308489	0.000642064842782305\\
};
\addlegendentry{$\mathcal{O}(h^2)$}

\end{axis}

\begin{axis}[%
width=3.486in,
height=4.245in,
at={(5.941in,0.573in)},
scale only axis,
xmode=log,
xmin=0.0462624701308489,
xmax=0.430042983911045,
xminorticks=true,
xlabel style={font=\color{white!15!black}},
xlabel={$h$},
ymode=log,
ymin=0.01,
ymax=2,
yminorticks=true,
axis background/.style={fill=white},
title style={font=\bfseries},
title={\large{Estimator and $L^2(H^1)$ error for different $\tau$}},
legend style={at={(0.01,0.99)}, anchor=north west, legend cell align=left, align=left, draw=white!15!black},
legend columns=2,
]
\addplot [color=white!30!black, mark=square, mark options={solid, white!30!black},line width=1pt,mark size=3pt]
  table[row sep=crcr]{%
0.430042983911045	0.993382769801442\\
0.348693200988337	0.413955658865301\\
0.180653453374279	0.195134914558749\\
0.0918101812671815	0.0767102213125241\\
0.0462624701308489	0.0551178112108131\\
};
\addlegendentry{$L^2(H^1)$ error $(\tau \!=\! 0.01)$}

\addplot [color=white!70!black, mark=square, mark options={solid, white!70!black},line width=1pt,mark size=3pt]
  table[row sep=crcr]{%
0.430042983911045	0.175695234977988\\
0.348693200988337	0.241583793078715\\
0.180653453374279	0.203741273948975\\
0.0918101812671815	0.167933967743723\\
0.0462624701308489	0.147768955099406\\
};
\addlegendentry{Estimator $\eta$ $(\tau \!=\! 0.01)$}

\addplot [color=white!30!black, mark=o, mark options={solid, white!30!black},line width=1pt,mark size=3pt]
  table[row sep=crcr]{%
0.430042983911045	0.704831555564995\\
0.348693200988337	0.157948758698773\\
0.180653453374279	0.0708927972455582\\
0.0918101812671815	0.0295899390446188\\
0.0462624701308489	0.0214491381919395\\
};
\addlegendentry{$L^2(H^1)$ error $(\tau \!=\! 0.001)$}

\addplot [color=white!70!black, mark=o, mark options={solid, white!70!black},line width=1pt,mark size=3pt]
  table[row sep=crcr]{%
0.430042983911045	0.207021323303769\\
0.348693200988337	0.19071651798131\\
0.180653453374279	0.114170704157295\\
0.0918101812671815	0.0659934381888048\\
0.0462624701308489	0.0411577757782331\\
};
\addlegendentry{Estimator $\eta$ $(\tau \!=\! 0.001)$}

\addplot [color=white!30!black, mark=x, mark options={solid, white!30!black},line width=1pt,mark size=3pt]
  table[row sep=crcr]{%
0.430042983911045	0.676090496953287\\
0.348693200988337	0.108040226554879\\
0.180653453374279	0.0450732033254962\\
0.0918101812671815	0.0232309133323613\\
0.0462624701308489	0.0187326388540756\\
};
\addlegendentry{$L^2(H^1)$ error $(\tau \!=\! 0.0001)$}

\addplot [color=white!70!black, mark=x, mark options={solid, white!70!black},line width=1pt,mark size=3pt]
  table[row sep=crcr]{%
0.430042983911045	0.214714895341424\\
0.348693200988337	0.20141392063202\\
0.180653453374279	0.115661015401895\\
0.0918101812671815	0.0610460721588877\\
0.0462624701308489	0.031337688478144\\
};
\addlegendentry{Estimator $\eta$ $(\tau \!=\! 0.0001)$}

\addplot [color=black, dashed]
  table[row sep=crcr]{%
0.430042983911045	0.430042983911045\\
0.348693200988337	0.348693200988337\\
0.180653453374279	0.180653453374279\\
0.0918101812671815	0.0918101812671815\\
0.0462624701308489	0.0462624701308489\\
};
\addlegendentry{$\mathcal{O}(h)$}

\end{axis}
\end{tikzpicture}%
%
    }
	\caption{Reliability and efficiency analysis for the droplet experiment.}
	\label{fig:rel_eff_plot}
\end{figure}

The error (in black) and the indicator (in grey) are presented in Figure~\ref{fig:rel_eff_plot} for the set of meshes and fixed time step sizes (differentiated by the different marks) against the spatial mesh-width $h$.

Since our estimators bound both error norms simultaneously with the same order, one can observe that for sufficiently small $\tau$ and $h$ the estimator follows the dominating $L^2(H^1)$-error curves. The efficiency is heuristically observed by noting that the error curves behave similar to the dominating error part. For small step sizes this is the $L^2(H^1)$-error, and for small temporal step sizes the $L^\infty(\mathcal{H}^{-1}_{L,\beta})$-error. 

Note that we scaled the $L^2(H^1)$-errors by $\eps$ in accordance to Theorem~\ref{thm:MainResultReliability}.

This indicates, while we did not prove it, that the derived error estimates are possibly not only reliable but also efficient. Thus they seem to mimic the error behavior accurately.

\subsection{Runtime analysis}
\label{sec:ex3}
In general, numerical computations for Cahn-Hilliard type problems tend to become increasingly difficult for decreasing $\eps$.
We cannot expect this to be overcome by the adaptive approach, however we want to showcase that adaptivity can decrease computational cost, which enlarges the set of $\eps$ which are feasible for computation.

The experiment setup is as follows:
For $\eps = 0.01$, $K = L = 1$, and $T = 0.001$ we computed a set of uniform and adaptive approximations, derived for different spatial meshes of the unit ball and different tolerances respectively. The initial data is given as an elliptically shaped droplet with barycenter at $(-1,0)$, a maximal horizontal elongation of $0.75$ and a maximal vertical elongation of $0.18$, as in \eqref{eq:Num initial data droplet}.
For the uniform approach we set a fixed temporal step size $\tau=10^{-6}$ which was chosen such that the numerical scheme does not diverge and a use meshes with degrees of freedom $20 \cdot 10^{k}$ for $k=0,\dots,9$. The adaptive solutions are obtained by our algorithm for a starting mesh with $160$ degrees of freedom, a maximal time step size of $\tau_{\max}=10^{-5}$, and tolerances
\begin{align*}
	\TOL_{h},\TOL_c &\in \{20, 10, 5, 2, 1, 0.7, 0.5\}, \\
	\TOL_{\tau} &\in \{1, 0.75, 0.625, 0.55, 0.525, 0.5175, 0.5125\}.
\end{align*}
All approximations are compared to a reference solution with a spatial resolution of $123457$ degrees of freedom and time step size of $2\cdot 10^{-7}$.

\begin{figure}[htbp]
	\centering
	\includegraphics[width=\textwidth]{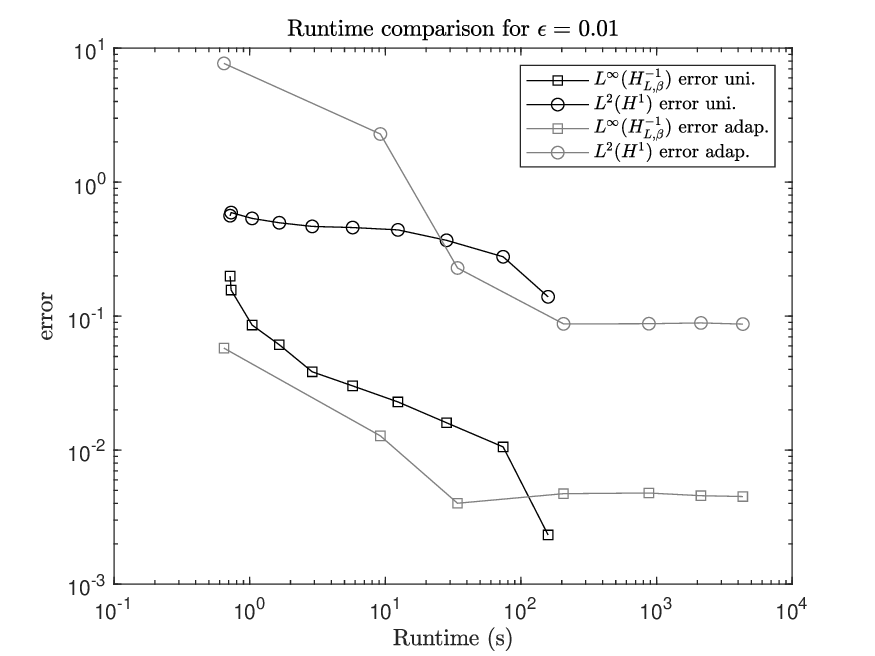}
	\caption{Runtime comparison between uniform and adaptive approximation of a reference solution.}
	\label{fig:runtimes}
\end{figure}

In Figure~\ref{fig:runtimes} we see the runtime plotted against the $L^\infty(H_{L,\beta}^{-1})$ (square marker) and $L^2(H^1)$ (circle marker) error for the uniform (black) and adaptive (grey) approximations. One can nicely see that our adaptive algorithm is in general more efficient, until the error becomes almost constant, which we believe to be the case because of the naturally limited accuracy of the reference solution. Additionally, we note that the apparent superconvergence of the uniform solution for the longest runtime arises due to the error-cancellation against the (structurally identically constructed) reference-solution.

\section*{Acknowledgments}
The authors would like to thank Stefan Metzger for insightful discussions on mass conservation properties of adaptive schemes, inspiring possible future research directions. We also thank Bal\'azs Kov\'acs for many valuable discussions and suggestions regarding the project. 

\section*{Declaration of AI-usage}

We declare that no artificial intelligence tools were used in the preparation of this document, except for the following two very limited cases: While preparing the runtime experiment, an AI tool (OpenAI: GPT5.3-Codex) was used to vectorize those parts of the algorithm, which compute the estimators. In the second case we used the AI tool in order to write a subprogram computing the error between the adaptive, respectively uniform, and the reference solution.

\appendix
\section{Estimates for interpolation bounds}
\label{sec:app_estimates_int_bounds}
We provide a detailed derivation of the two non-trivial bounds in \ref{Lemma:common_interpolation}.
I.e.\ we proof for $T \in \Omega_h^n$ that
\begin{align*}
\Big\|\frac{\d}{\d \theta} x^\theta - \big(\partial_\theta F^\theta\big)\big((F^\theta)^{-1}\big(x^\theta\big)\Big\|_{L^\infty(T^{\ell})} &\leq ch_{T}^2 \\
\Big\|D_a \Big(\frac{\d}{\d \theta} x^\theta- \big(\partial_\theta F^\theta\big)\big((F^\theta)^{-1}(x^\theta) \Big)\Big\|_{L^\infty(T^{\ell})} &\leq ch_{T}.
\end{align*}

\begin{proof}
First we note that each parent element $T \in \mathcal{T}_h^n$ has a set of corresponding elements $T_i^j \in \mathcal{T}_h^\theta$ in the intermediate mesh of $\Omega_h^\theta$.
We only focus on a single refinement/coarsening level in the following argument, since multiple levels are handled by a different argument (see \cite{KovacsLantelme2026} but represented on the parent mesh).
The parent element, effected by the lift, has two corresponding elements.
The following results are stated on the parent element however should be understood elementwise for each corresponding element.
This is done as the parent element dictates the bounds, since perturbation results including the lift are dependent on the element size of $T$.

Both estimates are bound by splitting the estimates and bounding each term separately.
\textbf{Proof of first estimate.}
We employ results for the lift map $G_h$ and its components $\Phi$ and $F$ introduced in \cite[Section~4]{ElliottRanner2013}. Given some small enough $h \leq h_0$ we use the following estimates with $h$-independent constant
\begin{equation}
\label{eq:ell_rann_bulk_estimates}
\begin{aligned}	
	\|G_h\|_{W^{2,\infty}(T)} &\leq c, \\
	\|(D_x G_h)^{-1}\|_{L^\infty(T)} &\leq c, \\
	\| D_{\hat{x}}^m \Phi \|_{L^\infty(\hat{T})} &\leq ch_{T}^2 \qquad \text{for} \quad m = 0,1,2.
\end{aligned}	
\end{equation}
These follow from \cite[Proposition~4.4, Lemma~4.6 \& Proposition~4.7]{ElliottRanner2013}, the inverse relation is obtained by a Neumann series argument.
	
We quantify $\frac{\d}{\d \theta} x^\theta(a)$ for $a \in T^{\ell}$ which measures the rate of change of an unlifted point of the bulk--surface with respect to the underlying change of the $\theta$ surface. 
	
	First note by differentiation of inverse functions that 
	\begin{equation}
	\label{eq:coarsening_proof_dtheta_inverse_G}
		\frac{\d}{\d \theta} x^\theta(a) = \frac{\d}{\d \theta} [(G_h^\theta)^{-1}(a)] = -\big((D_x G_h^\theta)(x^\theta(a))\big)^{-1} \cdot (\partial_\theta G_h^\theta)(x^\theta(a)),
	\end{equation}
	then we quantify the rate of change in the $\theta$ argument of the lift map, using $G_h^\theta = F^{e,\theta} \circ (F^\theta)^{-1}$  as:
	\begin{align*}
	(\partial_\theta G_h^\theta)(x) &= \big(D_{\hat{x}} F^{e,\theta}\big)\big((F^\theta)^{-1}(x)\big) \cdot \frac{\d}{\d \theta} (F^\theta)^{-1}(x) + (\partial_\theta F^{e,\theta})((F^\theta)^{-1}(x)) 
	\end{align*}

	We continue by expanding $F^{e,\theta} = F^\theta + \Phi^\theta$, explicitly computing the two partial derivatives and use the direct estimates \eqref{eq:ell_rann_bulk_estimates}. The situation and different contributions are depicted in Figure \ref{fig:theta-construction}. By the fact that only the node $x_l$ of an element changes and its rate of change is known to be $\mathcal{O}(h_{T}^2)$ the map 
	\[
	F^\theta = A(\theta) \hat{x} + b
	= \begin{bmatrix}
		\vert & \vert & \cdots & \vert & \cdots & \vert \\
		x_{N+1} - x_1 & x_N - x_1 & \cdots & x_l^\theta-x_1 &\cdots & x_2 - x_1 \\
		\vert & \vert & \cdots & \vert &  \cdots & \vert
	\end{bmatrix} \begin{bmatrix}
		\vert \\
		\hat{x} \\
		\vert
	\end{bmatrix} + 
	\begin{bmatrix}
		\vert \\
		x_1 \\
		\vert
	\end{bmatrix}
	\]
	is trivial to differentiate, thus for any point $\hat{x} \in \hat{T}$ we have
	\begin{equation}
		\label{eq:coarsening_proof_F_partialtheta}
		\begin{aligned}
		\frac{\d}{\d \theta} F^\theta(\hat{x}) &= \frac{\d}{\d\theta} A(\theta)\,\hat{x} = \mathcal{O}(h_{T}^2),\\
		D_{\hat{x}}F^\theta(\hat{x}) &= A(\theta) = \mathcal{O}(h_{T}).
		\end{aligned}
	\end{equation}
	This also allows us to bound the second term of \eqref{eq:coarsening_proof_thm_of_calculus} in the same manner, where the only difference is the evaluation point.

	Next we quantify the rate of change of $\Phi$ under the movement of the surface for a fixed reference node. The situation is depicted in the right side of Figure~\ref{fig:theta-construction}. We simplify the derivative along the second argument as follows 
	\begin{align}
	\label{eq:coarsening_proof_partial_dtheta_Phi}
		\frac{\d}{\d \theta} \Phi^\theta (\hat{x}) &= (\hat{\lambda}^\star)^{3} \partial_\theta (p(y^\theta)-y^\theta) = (\hat{\lambda}^\star)^{3} \big((D_{y} p)(y^\theta)-I \big) \partial_\theta y^\theta(\hat{x}).
	\end{align}
	Note that $D_{y} p(y)-I = D_y (p(y)-y) = \mathcal{O}(h_{T})$ from \cite{ElliottRanner2013}. 

	Further by $\partial_\theta y^\theta(\hat{x}) = \frac{\hat{\lambda}_l}{\hat{\lambda}^\ast} \frac{\d}{\d \theta} x_l^\theta = \mathcal{O}(h_{T}^2)$ and $\hat{\lambda}^\star \leq 1$,
	we are now able to also quantify the $\Phi$ derivatives. Based on using these bounds and \eqref{eq:ell_rann_bulk_estimates} for the spatial derivative we have
	\begin{align*}
	\frac{\d}{\d \theta} \Phi^\theta (\hat{x}) &= \mathcal{O}(h_{T}^3) \\
	D_{\hat{x}} \Phi^\theta (\hat{x}) &= \mathcal{O}(h_{T}^2).
	\end{align*}
	Finally by $A(\theta)^{-1} = \mathcal{O}(h_{T}^{-1})$ and $x-x_1 = \mathcal{O}(h_{T})$ we have for a fixed $x \in T$ the relation
	\begin{align*}
	\frac{\d}{\d\theta} (F^\theta)^{-1}(x) = \frac{\d}{\d\theta} A(\theta)^{-1} (x-x_1) = -A(\theta)^{-1} \Big(\frac{\d}{\d\theta} A(\theta) \Big) A(\theta)^{-1} (x-x_1) = \mathcal{O}(h_{T}).
	\end{align*}		
	In total these bounds show that 
	\begin{align}
		\label{eq:coarsening_proof_G_theta_derv}
		(\partial_\theta G_h^\theta)(x^\theta(a)) = \mathcal{O}(h_{T}^2).
	\end{align}

	The first estimate is the combination of the the derived bound \eqref{eq:coarsening_proof_F_partialtheta}, and by inserting \eqref{eq:ell_rann_bulk_estimates}, \eqref{eq:coarsening_proof_G_theta_derv} into \eqref{eq:coarsening_proof_dtheta_inverse_G}, and the assumption that $h \leq h_0$.	
	
	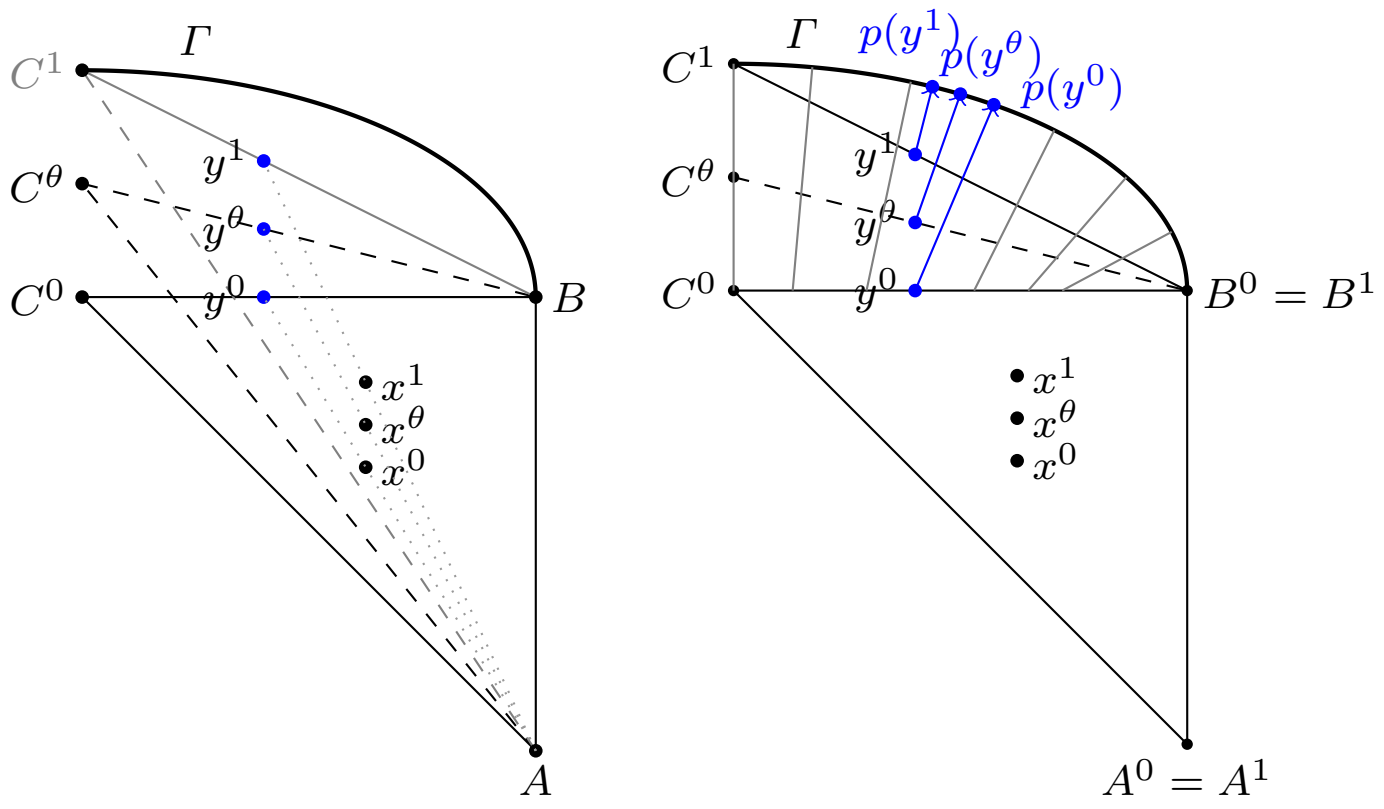
\begin{figure}[htbp]
\centering
\begin{tikzpicture}[scale=1.5]
  \draw[thick, smooth, samples=200, domain=0:90]
    plot({2*cos(\x)}, {sin(\x)});
  \node[above] at (0.5,1) {$\Gamma$};

  \coordinate (B)  at (2,0);
  \coordinate (A)  at (2,-2);
  \coordinate (C0) at (0,0);
  \coordinate (Cth) at (0,0.5);
  \coordinate (C1) at (0,1);

  \draw[black] (C0) node[left] {$C^0$} -- (B) node[right] {$B$};
  \draw[black,dashed] (Cth) node[left] {$C^\theta$} -- (B);
  \draw[gray] (C1) node[left] {$C^1$} -- (B);

  \draw[black] (A) node[below] {$A$} -- (B);
  \draw[black] (A) -- (C0);
  \draw[black,dashed] (A) -- (Cth);
  \draw[gray,dashed] (A) -- (C1);

  \fill (A)  circle[radius=0.03];
  \fill (B)  circle[radius=0.03];
  \fill (C0) circle[radius=0.03];
  \fill (Cth) circle[radius=0.03];
  \fill (C1) circle[radius=0.03];

	\coordinate (x0)  at ($ 3/8*(A) + 1/4*(B) + 3/8*(C0) $);
	\coordinate (xth) at ($ 3/8*(A) + 1/4*(B) + 3/8*(Cth) $);
	\coordinate (x1)  at ($ 3/8*(A) + 1/4*(B) + 3/8*(C1) $);
	
	\fill[black] (x0)  circle[radius=0.03];
	\fill[black] (xth) circle[radius=0.03];
	\fill[black] (x1)  circle[radius=0.03];
	
	\node[right] at (x0) {$x^0$};
	\node[right] at (xth) {$x^\theta$};
	\node[right] at (x1) {$x^1$};

  \path[name path=ray0]  (A) -- ($(A)!1.6!(x0)$);
  \path[name path=rayth] (A) -- ($(A)!1.6!(xth)$);
  \path[name path=ray1]  (A) -- ($(A)!1.6!(x1)$);

  \path[name path=edge0]  (B) -- (C0);
  \path[name path=edgeth] (B) -- (Cth);
  \path[name path=edge1]  (B) -- (C1);

  \path[name intersections={of=ray0 and edge0, by=y0}];
  \path[name intersections={of=rayth and edgeth, by=yth}];
  \path[name intersections={of=ray1 and edge1, by=y1}];

  \fill[blue] (y0)  circle[radius=0.03];
  \fill[blue] (yth) circle[radius=0.03];
  \fill[blue] (y1)  circle[radius=0.03];

  \node[left] at (y0) {$y^0$};
  \node[left]       at (yth) {$y^\theta$};
  \node[left] at (y1) {$y^1$};


  \draw[black!40,dotted] (A) -- ($(A)!1.6!(x0)$);
  \draw[black!40,dotted] (A) -- ($(A)!1.6!(xth)$);
  \draw[black!40,dotted] (A) -- ($(A)!1.6!(x1)$);

\end{tikzpicture}
~
		\begin{tikzpicture}[scale=1.5]
			\draw[thick, smooth, samples=200, domain=0:90]
			plot({2*cos(\x)}, {sin(\x)});
			\coordinate (top) at (0.3,1);
			\node[above] at (top) {$\Ga$};
			\draw[black] (0,0) node[left]{$C^0$} -- (2,0);
			\fill (0,0)    circle[radius=0.025];
			
			\draw[black] (0,1) node[left]{$C^1$} -- (2,0) node[right]{$B^0 = B^1$};
			\fill (2,0)    circle[radius=0.025];
			\fill (0,1)    circle[radius=0.025];
			
			\draw[black, dashed] (0,0.5) node[left]{$C^\theta$} -- (2,0);
			\fill (0,0.5)    circle[radius=0.025];		
			
			\coordinate (A) at (2,-2); 
			\draw[black] (A) node[below]{$A^0 = A^1$} -- (2,0);
			\draw[black] (A) -- (0,0);
			\fill (2,-2) circle[radius=0.025];		
			
			
			\path[name path=C0] (0,0) -- (2,0);
			\path[name path=C1] (0,1) -- (2,0);
			\path[name path=Ctheta] (0,0.5) -- (2,0);		
			
	\coordinate (x0)  at ($ 3/8*(A) + 1/4*(B) + 3/8*(C0) $);
	\coordinate (xth) at ($ 3/8*(A) + 1/4*(B) + 3/8*(Cth) $);
	\coordinate (x1)  at ($ 3/8*(A) + 1/4*(B) + 3/8*(C1) $);
	
	\fill[black] (x0)  circle[radius=0.03];
	\fill[black] (xth) circle[radius=0.03];
	\fill[black] (x1)  circle[radius=0.03];
	
	\node[right] at (x0) {$x^0$};
	\node[right] at (xth) {$x^\theta$};
	\node[right] at (x1) {$x^1$};

  \path[name path=ray0]  (A) -- ($(A)!1.6!(x0)$);
  \path[name path=rayth] (A) -- ($(A)!1.6!(xth)$);
  \path[name path=ray1]  (A) -- ($(A)!1.6!(x1)$);

  \path[name path=edge0]  (B) -- (C0);
  \path[name path=edgeth] (B) -- (Cth);
  \path[name path=edge1]  (B) -- (C1);

  \path[name intersections={of=ray0 and edge0, by=y0}];
  \path[name intersections={of=rayth and edgeth, by=yth}];
  \path[name intersections={of=ray1 and edge1, by=y1}];

  \fill[blue] (y0)  circle[radius=0.03];
  \fill[blue] (yth) circle[radius=0.03];
  \fill[blue] (y1)  circle[radius=0.03];

  \node[left] at (y0) {$y^0$};
  \node[left]       at (yth) {$y^\theta$};
  \node[left] at (y1) {$y^1$};
			
			\coordinate (py0) at ({2*cos(55)},{sin(55)});		
			\coordinate (pytheta) at ({2*cos(60)},{sin(60)});
			\coordinate (py1) at ({2*cos(64)},{sin(64)});
			\fill[blue] (py0) circle[radius=0.03];
			\fill[blue] (pytheta) circle[radius=0.03];
			\fill[blue] (py1) circle[radius=0.03];
			\node[above,blue,yshift=-4pt, xshift= 15pt] at (py0) {$p(y^0)$};
			\node[above,blue,yshift=0pt,xshift= 6pt] at (pytheta) {$p(y^\theta)$};
			\node[above,blue,yshift=3pt, xshift = -4pt] at (py1) {$p(y^1)$};
			\draw [blue, ->] (y0) -- (py0);
			\draw [blue, ->] (yth) -- (pytheta);
			\draw [blue, ->] (y1) -- (py1);
			
			\foreach \t in {15,30,45,67,80,90} {
				\pgfmathsetmacro{\xi}{2*cos(\t)}
				\pgfmathsetmacro{\yi}{sin(\t)}
				\pgfmathsetmacro{\tx}{-2*sin(\t)}
				\pgfmathsetmacro{\ty}{ cos(\t)}
				\pgfmathsetmacro{\nx}{\ty}
				\pgfmathsetmacro{\ny}{-\tx}
				\pgfmathsetmacro{\szero}{-\yi/\ny}
				\pgfmathsetmacro{\xzero}{\xi + \szero*\nx}
				\draw[gray] (\xi,\yi) -- (\xzero, 0);
				\ifnum\t=45
				\pgfmathsetmacro{\sone}{(2 - \xi - 2*\yi)/(\nx + 2*\ny)}
				\pgfmathsetmacro{\uone}{(\xi + \sone*\nx)/2}
				\pgfmathsetmacro{\xone}{2*\uone}
				\pgfmathsetmacro{\yone}{1 - \uone}

				
				\fi			
			}
		\end{tikzpicture}
\caption{Sketch of the two-dimensional situation. For the element $\Delta ABC^\theta$ of $\Omega_h^\theta$ with $\theta\in[0,1]$, we fix $\hat x = (3/8,1/4, 3/8)$. The left image shows the linear change in $x^\theta = F^\theta(\hat{x})$ and $y^\theta$ with respect to $\theta$. The right figure shows the non-linear change in $p(y^\theta)$ defined by the unique point projection.}
\label{fig:theta-construction}
\end{figure}	

\textbf{Proof of second estimate.}	
We directly compute
	\begin{align}
	\label{eq:coarsening_proof_stiffness_Da_F_inv}
	D_a \Big( \big(\partial_\theta F^\theta\big)\big((F^\theta)^{-1}(x^\theta(a))\big) \Big) &= D_a \Big( \frac{\d}{\d \theta} A(\theta) \cdot (F^\theta)^{-1}(x^\theta(a))\Big) \nonumber \\
	&= \frac{\d}{\d \theta} A(\theta) \cdot \big(D_x (F^\theta)^{-1}\big)\big(x^\theta(a)\big) \cdot D_a x^\theta(a) \nonumber  \\
	&= \frac{\d}{\d \theta} A(\theta) A(\theta)^{-1} [(D_x G_h^\theta)(D_a x^\theta(a))]^{-1} =  \mathcal{O}(h_{T}),
	\end{align}	 
	which employs the explicit description of $F^\theta$, the bounds to the reference matrix $A$, namely $|A^{-1}| = \frac{c}{h_{T}}$ elementwise  \cite{ElliottRanner2013} and $A'(\theta) = \mathcal{O}(h_{T}^2)$, and the derivative bound of the lift map $G_h^\theta$ by \eqref{eq:ell_rann_bulk_estimates}. 

	The mixed derivative of the lift map needs further analysis. We insert the formula \eqref{eq:coarsening_proof_dtheta_inverse_G} and apply product and chain rule to obtain
	\begin{align}
	\label{eq:coarsening_proof_mixed_deriv_G}
	  	D_a \Big(\frac{\d}{\d \theta} x^\theta(a)\Big) = -\Bigg(&D_x \Big[D_x (G_h^\theta)^{-1}\Big] (x^\theta(a)) \cdot \big(D_a x^\theta\big)(a) \cdot (\partial_\theta G_h^\theta)(x^\theta(a)) \nonumber \\
		&+ \big((D_x G_h^\theta)(x^\theta(a))\big)^{-1} \cdot \big(D_{x}(\partial_\theta G_h^\theta)\big)(x^\theta(a)))\big) \cdot \big(D_a x^\theta\big)(a) \Bigg),
	\end{align}
	where we again identified $x^\theta(a) = (G_h^\theta)^{-1}(a)$ in evaluation points.
	
	We go over the orders for each of the terms, the second order spatial derivative is simplified using the component wise differentiation of matrix inverses and by \eqref{eq:ell_rann_bulk_estimates} to show
	\begin{align*}
		\|\partial_{x_j} \big[\big(D_x G_h^\theta\big)^{-1}\big]\| = \|-(D_x G_h^\theta)^{-1} (\partial_{x_j} D_x G_h^\theta) (D_x G_h^\theta)^{-1}\| \leq c.
	\end{align*}
	We also note that by \eqref{eq:ell_rann_bulk_estimates} we have both $\|(D_x G_h^\theta)^{-1}\| \leq 1$ and as an immediate consequence by inverse differentiation
	\begin{align*}	
		 \|D_a x^\theta\| = \|D_a (G_h^\theta)^{-1}\| \leq c.
	\end{align*}
	From the analysis of the mass terms we additionally derive by \eqref{eq:coarsening_proof_G_theta_derv} that $\partial_\theta G_h^\theta(x^\theta(a)) = \mathcal{O}(h_{T}^2)$.

The final piece is the bound for the mixed derivative of the lift map $G_h^\theta(x)$ we can simplify
	\begin{align*}
	\big(D_{x}(\partial_\theta G_h^\theta)\big)(x^\theta(a)) =  D_x \Big[ \big(D_{\hat{x}} F^{e,\theta}\big)\big((F^\theta)^{-1}(x)\big) \cdot \big(\partial_\theta (F^\theta)^{-1}\big)(x) + \big(\partial_\theta F^{e,\theta}\big)\big((F^\theta)^{-1}(x)\big) \Big] \Big|_{x =x^\theta(a)}.
	\end{align*}
	Next we can expand $F^{e,\theta} = F^\theta + \Phi^\theta$ and group the terms, while ignoring the evaluation point for now, to derive 
	\begin{align*}
	\big(D_{x} (\partial_\theta G_h^\theta)\big)(x) &=  \big(D^2_{\hat{x}} (F^\theta + \Phi^\theta)\big)\big((F^\theta)^{-1}(x)\big) \cdot \big(D_x (F^\theta)^{-1}\big)(x) \cdot \big(\partial_\theta (F^\theta)^{-1}\big)(x) \\
	&\quad +  \big(D_{\hat{x}} (F^\theta + \Phi^\theta)\big)\big((F^\theta)^{-1}(x)\big) \cdot (D_x (\partial_\theta (F^\theta)^{-1}))(x)	\\
	&\quad + (D_{\hat{x}} \partial_\theta (F^\theta + \Phi^\theta))((F^\theta)^{-1}(x)) \cdot \big(D_x (F^\theta)^{-1}\big)(x).  
	\end{align*}
	Based on the fact that the affine transformation $F^\theta(x) = A(\theta)x + b$ and its inverse are easily differentiable we can simplify
	\begin{align*}
	D_{x^\theta}\big((\partial_\theta G_h^\theta)(x^\theta)\big) &=  \big(D^2_{\hat{x}} \Phi^\theta\big)\big((F^\theta)^{-1}(x^\theta)\big) \cdot A(\theta)^{-1} \Big(\frac{\d}{\d \theta} A(\theta)^{-1}\Big) (x^\theta-b)  \\
	&\quad +  \Big(A(\theta) + \big(D_{\hat{x}} \Phi^\theta\big)\big((F^\theta)^{-1}(x^\theta)\big)\Big) \cdot \Big(\frac{\d}{\d \theta} A(\theta)^{-1}\Big)\\
	&\quad + \Big(D_{\hat{x}} \partial_\theta (F^\theta + \Phi^\theta)\Big)((F^\theta)^{-1}(x^\theta)) \cdot A(\theta)^{-1}.  
	\end{align*}
	
	We have already quantified the size of the quantities related to the affine map, additionally by matrix inverse differentiation we have $|\frac{\d}{\d \theta}A(\theta)^{-1}| \leq c$.
	Further we can employ the bounds of \eqref{eq:ell_rann_bulk_estimates} to bound the spatial derivatives of $\Phi^\theta$.
	To bound the full expression we have to determine the order of the mixed derivatives.
	The differentiation of the affine map simply yields 
	$$(D_{\hat{x}} \partial_\theta F^\theta)(\hat{x}) = \frac{\d}{\d \theta} A(\theta) = \mathcal{O}(h_{T}^2).$$
	
	Next, the explicit differentiation of $\Phi$  by using \eqref{eq:coarsening_proof_partial_dtheta_Phi}:
	\begin{align}
	\label{eq:coarsening_proof_mixed_deriv_Phi}
	(\partial_\theta D_{\hat{x}} \Phi^\theta)(\hat{x}) = D_{\hat{x}} \Big((\hat{\lambda}^\star)^3 \big((D_{y} p)(y^\theta)-I\big) \partial_\theta y^\theta(\hat{x}) \Big)
	\end{align}
	We have already establishes bounds of each component, without the additional spatial differentiation.
	Determining and bounding each component separately lets use bound the full derivative. We have, due the linearity of the barycentric coordinates $\hat{\lambda}_i$
	\begin{align*}
		D_{\hat{x}} (\hat{\lambda}^\star)^3 = 3 (\hat{\lambda}^\star)^2 \sum_{i=1}^l D_{\hat{x}} \hat{\lambda}_i \leq c(\hat{\lambda}^\star)^2.
	\end{align*}
	By quotient rule, that the vertices $x_i$ are fixed for all $i \neq l$, and the linearity and boundedness of $\hat{\lambda}_i$ we also have
	\begin{align*}
		D_{\hat{x}} \partial_\theta y^\theta(\hat{x}) = D_{\hat{x}}\Big(\frac{\hat{\lambda}_l}{\hat{\lambda}^\star}\Big) \frac{\d}{\d \theta} x_l(\theta) = \frac{(D_{\hat{x}} \hat{\lambda}_l) \hat{\lambda}^\star - \hat{\lambda}_l D_{\hat{x}} \hat{\lambda}^\star}{(\hat{\lambda}^\star)^2} \frac{\d}{\d \theta} x_l(\theta) \leq c h_{T}^2 \frac{1}{{(\hat{\lambda}^\star)^2}}.
	\end{align*}
	Finally, by $\|D_y^2 p(y)\|_{L^\infty} = \|D_y^2 \big(p(y)-y\big)\|_{L^\infty} \leq c$ (within the proof of \cite[Lemma~4.3]{ElliottRanner2013}) and $\|D_{\hat{x}} y\|_{L^\infty} \leq c \frac{h}{\lambda^\star}$ (see \cite[Lemma~6.3]{Bernardi1989}) we bound
	\begin{align*}
		D_{\hat{x}}((D_{y} p)(y^\theta)-I) = (D_{y}^2 p)(y^\theta) \cdot (D_{\hat{x}} y^\theta)(\hat{x})  \leq c \frac{h}{\lambda^\star}. 
	\end{align*}
	Now we can apply product rule to \eqref{eq:coarsening_proof_mixed_deriv_Phi} and insert our bounds to get
	\begin{align*}
		(\partial_\theta D_{\hat{x}} \Phi^\theta)(\hat{x}) = \mathcal{O} (h_{T}^3).
\end{align*}	 	
	We observe that the problematic quotients of powers $\hat{\lambda}^\star$ cancel when combining the separate bounds which guarantees that the estimates do not degenerate for points close to the non-boundary vertices. 
	
	Now this implies in total that $\|D_x \partial_\theta G_h^\theta\| \leq c h_{T}$ and comparing to \eqref{eq:coarsening_proof_mixed_deriv_G} shows that the lift-related contribution to the error for the second estimate is of first order
	\begin{align}
		\label{eq:coarsening_proof_stiffness_bound_Da_d_partial_theta_G_inv}
		\Big\|D_a \Big(\frac{\d}{\d \theta} x^\theta(a)\Big)\Big\| \leq ch_T.
	\end{align}

	Thus we have shown for both terms of the second estimate that they are bound by first order contributions, see \eqref{eq:coarsening_proof_stiffness_Da_F_inv} and \eqref{eq:coarsening_proof_stiffness_bound_Da_d_partial_theta_G_inv}, which finishes the proof. 
\end{proof}

\bibliographystyle{alpha}
\bibliography{C-Heq_adap}

\end{document}